\documentclass[12pt]{amsart}
\usepackage{amscd}
\usepackage{verbatim}
\usepackage{amssymb}

\newtheorem{thm}{Theorem}[section]
\newtheorem{prop}[thm]{Proposition}

\newtheorem{lemma}[thm]{Lemma}
\newtheorem{cor}[thm]{Corollary}
\newtheorem{con}[thm]{Conjecture}

\theoremstyle{definition}
\newtheorem{rem}[thm]{Remark}
\newtheorem{example}[thm]{Example}
\newtheorem{dfn}[thm]{Definition}

\newcommand{\onto}{\rightarrow\!\!\rightarrow}

\renewcommand{\Im}{\mathop{\mathrm{Im}}}

\newcommand{\CH}{\mathop{\mathrm{Ch}}\nolimits}

\newcommand{\BCH}{\mathop{\bar{\mathrm{Ch}}}\nolimits}
\newcommand{\BCHE}{\mathop{\mathrm{Che}}\nolimits}

\newcommand{\res}{\mathop{\mathrm{res}}\nolimits}

\newcommand{\pr}{\operatorname{\mathit{pr}}}
\newcommand{\inc}{\operatorname{\mathit{in}}}

\newcommand{\Char}{\mathop{\mathrm{char}}\nolimits}

\newcommand{\Z}{\mathbb{Z}}

\newcommand{\Spec}{\operatorname{Spec}}

\newcommand{\an}{0}

\newcommand{\Sum}{\operatornamewithlimits{\textstyle\sum}}
\newcommand{\Oplus}{\operatornamewithlimits{\textstyle\bigoplus}}

\newcommand{\<}{\left<}
\renewcommand{\>}{\right>}
\renewcommand{\ll}{\<\!\<}
\newcommand{\rr}{\>\!\>}

\newcommand{\compose}{\circ}

\newcommand{\IW}{\mathfrak{i}_0}
\newcommand{\iw}{\mathfrak{i}}
\newcommand{\jw}{\mathfrak{j}}

\newcommand{\Hight}{\mathfrak{h}}

\newcommand{\SP}{\operatorname{sp}}

\newcommand{\Sym}{\operatorname{Sym}}
\newcommand{\D}{D}

\newcommand{\BC}{\ast}
\newcommand{\NBC}{\compose}

\newcommand{\NC}{\compose}
\newcommand{\AC}{\bullet}

\title
{Holes in $I^n$}

\keywords
{Quadratic forms, Witt indices,
Chow groups, Steenrod operations, correspondences.
{\em 2000 Mathematical Subject Classifications:}
11E04; 14C25}

\date{December 2003}

\author
[N. Karpenko]
{Nikita A. Karpenko}

\address
{Laboratoire des Math{\'e}matiques\\
Facult{\'e} des Sciences Jean Perrin\\
Universit{\'e} d'Artois\\
rue Jean Souvraz SP 18\\
62307 Lens Cedex\\
France}

\email
{karpenko@euler.univ-artois.fr}

\begin{document}

\begin{abstract}
Let $F$ be an arbitrary field of characteristic $\ne2$.
We write $W(F)$ for the Witt ring of $F$ (for the definition, see
\cite[def. 1.9 of ch. 2]{Scharlau}), consisting of the isomorphism
classes of all anisotropic quadratic forms over $F$.
For any $x\in W(F)$, {\em dimension} $\dim x$ is defined as the
dimension of a quadratic form representing $x$.
The elements of all even dimensions form an ideal denoted $I(F)$.
The filtration of the ring $W(F)$ by the powers $I(F)^n$ of this ideal
plays a fundamental role in the algebraic theory of quadratic forms.
The Milnor conjectures, recently proved by Voevodsky and
Orlov-Vishik-Voevodsky, describe the successive quotients
$I(F)^n/I(F)^{n+1}$ of this filtration, identifying them with Galois cohomology
groups and with the Milnor $K$-groups modulo $2$ of the field $F$.
In the present article we give a complete answer to a different
old-standing question concerning $I(F)^n$, asking about the
possible values of $\dim x$ for $x\in I(F)^n$.
More precisely, for any $n\geq1$, we prove
that
$$
\dim I^n=\{2^{n+1}-2^i|\;i\in[1,\;n+1]\}\cup
\big(2\Z\cap[2^{n+1},\;+\infty)\big)\;,
$$
where $\dim I^n$ is the set of all $\dim x$ for all $x\in I(F)^n$ and
{\em all} $F$.
Previously available partial informations on $\dim I^n$ include the
classical Arason-Pfister theorem
\cite[Hauptsatz]{Arason-Pfister}
(saying that $(0,\;2^n)\cap\dim
I^n=\emptyset$) as well as a recent Vishik's theorem
\cite[th. 6.4]{Vishik-Lens}
on $(2^n,\;2^n+2^{n-1})\cap\dim I^n=\emptyset$
(the case $n=3$ is due to Pfister, \cite[Satz 14]{Pfister};
$n=4$ to Hoffmann, \cite[main th.]{Hoffmann-I^4}).
Our proof is based on computations in Chow groups of
powers of projective quadrics (involving the Steenrod
operations); the method developed can be also applied to other
types of algebraic varieties.
\end{abstract}

\maketitle

\tableofcontents

\section
{Introduction}
\label{Introduction}

In this text $F$ is a field with $\Char(F)\ne2$.
Let $I^n$ for some $n\geq1$
be the $n$-th power of the fundamental ideal $I$
(of the classes of the even-dimensional quadratic forms)
of the Witt ring $W(F)$.
An old-standing question in the algebraic theory of quadratic forms
asks about the possible values of
dimension of an anisotropic quadratic form $\phi$ over $F$ such that
$[\phi]\in I^n$, where $[\phi]$ is the class of $\phi$ in $W(F)$.

Examples with $\dim(\phi)=2^{n+1}-2^i$ for each $i\in\{1,2,\dots,n+1\}$
are easy to construct
(see Remark \ref{easy examples}).
A classical theorem of J. Arason and A. Pfister
\cite[Hauptsatz]{Arason-Pfister}
states that $\dim(\phi)$ is never between
$0=2^{n+1}-2^{n+1}$ and $2^n=2^{n+1}-2^n$.
Also it is known that every value between $2^n$ and $2^n+2^{n-1}=2^{n+1}-2^{n-1}$
is impossible (A. Pfister for $n=3$, \cite[Satz 14]{Pfister};
D. Hoffmann for $n=4$, \cite[main thm.]{Hoffmann-I^4};
A. Vishik for all $n$, \cite{Vishik-Lens}, see also
\cite[thm. 1.5]{Vishik-In} \footnote{Alternative proofs are given in
\cite[thm. 4.4]{KM} and in \cite{third proof}.}).

Finally, A. Vishik has shown that all even values $\geq2^{n+1}$ are possible
(\cite[thm. 4.12]{Vishik-In}, see also \S\ref{Possible dimensions} here)
and suggested the following

\begin{con}[{Vishik, \cite[conject. 4.11]{Vishik-In}}]\footnote{A.
Vishik announced (for the first time in June 2002)
that he has a proof
of Conjecture \ref{conj}; his proof however is not available.}
\label{conj}
If $\phi\in I^n$ and $\dim(\phi)< 2^{n+1}$, then $\dim(\phi)=2^{n+1}-2^i$
for some $i\in\{1,2,\dots,n+1\}$.
\end{con}

In the present text we prove this conjecture
(see \S\ref{Proof of Conjecture}),
obtaining a complete answer
to the question about possible dimensions of anisotropic quadratic forms
whose classes lie in $I^n$.
The proof closely follows the method of \cite{third proof}, but involves essentially
more computations.
As \cite{third proof} as well, it makes use of an important property of the quadratic forms
satisfying the hypotheses of Conjecture \ref{conj}  established
by A. Vishik in \cite{Vishik-In}.
Here we give an extended version of this result (see Proposition \ref{known})
with an elementary, complete, and self-contained (in particular,
independent of \cite{Vishik-In})
proof.

In the proof of Conjecture \ref{conj} we work with projective
quadrics rather than with quadratic forms themselves.
The method of proof is explained in \S \ref{Strategy of proof};
it certainly applies to other types of algebraic varieties
(in place of quadrics).

\bigskip
\noindent
{\sl Acknowledgements.}
This work was completed during the last week of September 2003
when I was visiting the University of Alberta;
an important part of this work was written down during the first week
of November 2003 when I was visiting the Universit{\"a}t G{\"o}ttingen,
where I gave a course of lectures on the subject.
I am very grateful to both institutions for hospitality and stimulating
conditions.
Besides I would like to mention that numerous remarks of A. S. Merkurjev
improved very much the final version of this text.

\section
{Notation and preliminary observations}
\label
{Notation and preliminary observations}

Everywhere in the text, $X$ is a smooth projective quadric over $F$ of an even dimension
$D=2d$ or of an odd dimension
\footnote{Only the even-dimensional case is important for our main purpose; the odd-dimensional
case is included for the sake of completeness.}
$D=2d+1$
given by a non-degenerate quadratic form $\phi$.
We write $X^r$ for the direct
product $X\times\dots\times X$ (over $F$) of $r$ copies of $X$
and we write $\BCH(X^r)$ for the image of the restriction homomorphism
$\CH(X^r)\to\CH(\bar{X}^r)$ where $\bar{X}=X_{\bar{F}}$ with a fixed algebraic closure
$\bar{F}$ of $F$ and $\CH(.)$ stays for the {\em modulo 2} total Chow group.
We say that an element of $\CH(\bar{X}^r)$ is {\em rational}, if it lies in the subgroup
$\BCH(X^r)\subset\CH(\bar{X}^r)$.

A basis of the group $\CH(\bar{X})$ (over $\Z/2\Z$) consists of
$h^i$ and $l_i$, $i\in\{0,1,\dots,d\}$, where $h$ stays for the
class of a hyperplane section of $\bar{X}$ while $l_i$ is the
class of an $i$-dimensional linear subspace\footnote{In the case
of even $D$ and $i=d$ (and only in this case) the class $l_i$
depends on the choice of the subspace: more precisely, there are
two different classes of $d$-dimensional subspaces on $\bar{X}$
and no canonical choice of one of them is possible; we do not care
about this however and we just choose one of them, call it $l_d$
and ``forget'' about the other one which is equal to $h^d-l_d$.}
lying on $\bar{X}$. Moreover, a basis of the group
$\CH(\bar{X}^r)$ for every $r\geq 1$ is given by the external
products of the basis elements of $\CH(\bar{X})$ (see, e.g.,
\cite[\S1]{i_1} for an explanation why this is a basis). Speaking
about a {\em basis} or a {\em basis element} (or a {\em basis
cycle}: we will often apply the word ``cycle'' to an element of a
Chow group) of $\CH(\bar{X}^r)$, we will always refer to the basis
described above. By the {\em decomposition} of an element
$\alpha\in\CH(\bar{X}^r)$ we always mean its representation as a
sum of basis cycles. We say that a basis cycle $\beta$ is {\em
contained} in the decomposition of $\alpha$ (or simply ``is
contained in $\alpha$''), if $\beta$ is a summand of the
decomposition. More generally, for two cycles
$\alpha',\alpha\in\CH(\bar{X}^r)$, we say that $\alpha'$ {\em is
contained} in $\alpha$, or that $\alpha'$ is a {\em subcycle} of
$\alpha$ (notation: $\alpha'\subset\alpha$), if every basis
element contained in $\alpha'$ is also contained in $\alpha$.

A basis element of $\CH(\bar{X}^r)$ is called {\em non-essential},
if it is an external product of powers of $h$ (including
$h^0=1=[\bar{X}]$); the other basis elements are called {\em
essential}. An element of $\CH(\bar{X}^r)$ which is a sum of
non-essential basis elements, is called non-essential as well.
Note that all non-essential elements are rational simply because
$h$ is rational.\footnote{There are at least two direct ways to
show that $h$ is rational: (1) $h$ is the pull-back of the
hyperplane class $H$ with respect to the embedding of $\bar{X}$
into the projective space, and $H$ is rational; (2) $h$ is the
first Chern class of $[\mathcal{O}_{\bar{X}}(1)]\in K_0(\bar{X})$,
and $[\mathcal{O}_{\bar{X}}(1)]=\res([\mathcal{O}_X(1)])$ is
rational.}

The multiplication table for the ring $\CH(\bar{X})$
is determined by the rules
$h^{d+1}=0$, $h\cdot l_i=l_{i-1}$ ($i\in\{d,d-1,\dots,0\}$;
we adopt the agreement that
$l_{-1}=0$), and
$l_d^2=(D+1)\cdot(d+1)\cdot l_0$ (see \cite[\S2.1]{alg.geom.inv.}
and \cite[\S 1.2.1]{disser}).
The multiplication tables for the rings $\CH(\bar{X}^r)$ (for all $r\geq2$)
follows by
$$
(\beta_1\times\dots\times\beta_r)\cdot
(\beta'_1\times\dots\times\beta'_r)=
\beta_1\beta'_1\times\dots\times\beta_r\beta'_r\;.
$$

The cohomological action of the topological Steenrod algebra on $\CH(\bar{X}^r)$
(see \cite{Brosnan} for the construction of the action of the topological Steenrod algebra
on the Chow group of a smooth projective variety;
originally Steenrod operations in algebraic geometry
were introduced (in the wider context of motivic
cohomology)
by V. Voevodsky, \cite{Voevodsky-operations})
is determined by
the fact that the total Steenrod operation
$S\!:\CH(\bar{X}^r)\to\CH(\bar{X}^r)$ is a
(non-homogeneous) ring homomorphism,
commuting with the external products
and satisfying the formulae (see \cite[\S2 and cor. 3.3]{i_1})
$$
S(h^i)=h^i\cdot(1+h)^i,\;\;
S(l_i)=l_i\cdot(1+h)^{D-i+1},\;\;
i\in\{0,1,\dots,d\}
$$
(in order to apply these formulae, one needs a computation of the binomial coefficients
modulo $2$, done, e.g., in \cite[lemma 1.1]{KMD}).

The group $\BCH(X)$ is easy to compute.
First of all one has

\begin{lemma}
\label{springer}
If the quadric $X$ is anisotropic
(that is, $X(F)=\emptyset$), then $\BCH_0(X)\not\ni l_0$.
\end{lemma}

\begin{proof}
If $l_0\in\BCH_0(X)$, then the variety $X$ contains a closed point
$x$ of an odd degree $[F(x):F]$.
It follows that the quadratic form $\phi$ is isotropic over
an odd degree extension of the base field (namely, over $F(x)$) and
therefore, by the Springer-Satz (see \cite[thm. 5.3]{Scharlau}),
is isotropic already over $F$.
\end{proof}

\begin{cor}
\label{X is non-essential}
If $X$ is anisotropic, then the group $\BCH(X)$ is generated by
the non-essential basis elements.
\end{cor}

\begin{proof}
If the decomposition of an element $\alpha\in\BCH(X)$ contains an
essential basis element $l_i$ for some $i\ne D/2$, then $l_i\in\BCH(X)$
because $l_i$ is the $i$-dimensional homogeneous component of $\alpha$
(and $\BCH(X)$ is a graded subring of $\CH(\bar{X})$).
If the decomposition of an element $\alpha\in\BCH(X)$ contains
the
essential basis element $l_i$ for $i= D/2$, then
the ($D/2$)-dimensional homogeneous component of $\alpha$
is either $l_{D/2}$ or $l_{D/2}+h^{D/2}$ and we still have $l_i\in\BCH(X)$.
It follows that $l_0=l_i\cdot h^i\in\BCH(X)$, a contradiction with
Lemma \ref{springer}.
\end{proof}

Now assume for a moment that the quadric $X$ is isotropic but not
completely split
(that is, $\IW(X)\leq d$),
write $a$ for the Witt index $\IW(X)$ of $X$
(defined as the Witt index $\IW(\phi)$ of $\phi$, see \cite[def. 5.10 of ch. 1]{Scharlau}),
and let $X_0$ be the projective
quadric
given by the anisotropic part $\phi_0$ of $\phi$
(one has $\dim(X_0)=\dim(X)-2a$; the case $\dim(X_0)=0$ is possible).
We consider a group
homomorphism $pr\!:\CH(\bar{X})\to\CH(\bar{X}_0)$ determined on the basis
by the formulae
$h^i\mapsto h^{i-a}$  and $l_i\mapsto l_{i-a}$
(here we adopt the agreement $h^i=0$ and $l_i=0$ for all negative $i$).
Also we consider a backward group homomorphism
$\inc\!:\CH(\bar{X}_0)\to\CH(\bar{X})$ determined by the formulae
$h^i\mapsto h^{i+a}$ and $l_i\mapsto l_{i+a}$ for
$i\in\{0,1,\dots,d-a\}$.

Let $r$ be a positive integer.
For every length $r$ sequence $i_1,\dots,i_r$ of integers satisfying
$i_j\in[0,\;a]\cup[D-a+1,\;D]$,
we define a group homomorphism
$$
\pr_{i_1\dots i_r}\!:\CH(\bar{X}^r)\to\CH(\bar{X}_0^s)
$$
with $s=\#\{i_j|\;i_j=a\}$,
called {\em projection}.
Let $\{j_1<\dots<j_s\}$ be the set of indices such that
$i_{j_s}=a$. We put
$J_l=\{j|\;i_j<a\}$ and
$J_h=\{j|\;i_j>a\}$.
Then we define
$\pr_{i_1\dots i_r}(\alpha_1\times\dots\times\alpha_r)$
for a basis element
$\alpha_1\times\dots\times\alpha_r$ as
$\pr(\alpha_{j_1})\times\dots\times\pr(\alpha_{j_s})$
as far as
$\alpha_j=l_{i_j}$ for any $j\in J_l$ and
$\alpha_j=h^{D-i_j}$ for any $j\in J_h$;
we set
$\pr_{i_1\dots i_r}(\alpha_1\times\dots\times\alpha_r)=0$
otherwise.

Also we define a backward group homomorphism
$\inc_{i_1\dots i_r}\!:\CH(\bar{X}_0^s)\to\CH(\bar{X}^r)$,
called {\em inclusion}, by
$$
\inc_{i_1\dots i_r}(\beta_1\times\dots\times\beta_s)=
\alpha_1\times\dots\times\alpha_r
$$
for a basis element $\beta_1\times\dots\times\beta_s$,
where $\alpha_j=l_{i_j}$ for $j\in J_l$,
$\alpha_j=h^{D-i_j}$ for $j\in J_h$, and
$\alpha_{j_k}=\inc(\beta_k)$ for $k=1,2,\dots,s$.

\begin{prop}[{cf. \cite[lemma 2.2]{third proof}}]
\label{isotropic}
In the notation introduced right above,
the direct sum
$$
\Oplus_{(i_1,\dots,i_r)}\pr_{i_1\dots i_r}\!:\CH(\bar{X}^r)\to\Oplus_{(i_1,\dots,i_r)}
\CH(\bar{X}_0^{s(i_1,\dots,i_r)})
$$
of all projections is an isomorphism
with the inverse given by the sum of all inclusions.
Under these isomorphisms, rational cycles correspond to rational
cycles.
\end{prop}

\begin{proof}
By the Rost motivic decomposition theorem for isotropic quadrics
(original proof is in \cite{Rost}, generalizations are obtained in
\cite{flag} and \cite{Chernousov-Gille-Merkurjev}),
there is a motivic decomposition (in the category of the integral Chow motives)
\begin{equation}
\tag{$*$}
X\simeq\Z\oplus\Z(1)\oplus\dots\oplus\Z(a-1)\oplus X_0(a)
\oplus\Z(D-a+1)\oplus\dots\oplus\Z(D)
\end{equation}
(where $\Z$ is the motive of $\Spec F$, while $M(i)$ is the $i$-th Tate twist of
a motive $M$).
Rasing to the $r$-th power, we get a motivic decomposition of the variety
$X^r$;
each summand of this decomposition is a twist of the motive of $X_0^s$ with $s$
varying between $0$ and $r$.
If we numerate the summands of the decomposition $(*)$
by their twists,
then the summands of the decomposition of $X^r$ are numerated by the sequences
$$
i_1,\dots,i_r\;\;\text{ with }\;i_j\in[0,\;a]\cup[D-a+1,\;D]\;.
$$
Moreover, the $(i_1,\dots,i_r)$-th summand is $X_0^s(i_1+\dots+i_r)$, where
$s=\#\{i_j|\;i_j=a\}$.

In order to finish the proof of Proposition \ref{isotropic}, it suffices to show that the
projection morphism to the $(i_1,\dots,i_r)$-th summand considered on the Chow group and over
$\bar{F}$ coincides with $\pr_{i_1\dots i_r}$ while the inclusion morphism of
the $(i_1,\dots,i_r)$-th summand considered on the Chow group and over
$\bar{F}$ coincides with $\inc_{i_1\dots i_r}$.
Clearly it suffices to check this for $r=1$ only.
For $i\ne D/2$,
this is particularly easy to do because of the relation
$\dim_{\Z/2\Z}(\CH_i(\bar{X}))\leq 1$.
Indeed, $\CH_i(\Z(k))=0$ for $k\ne i$.
Therefore for any $i$ with $a\leq i\leq D-a$, $i\ne D/2$,
the projection and the inclusion between $\CH_i(\bar{X})$ and $\CH_{i-a}(\bar{X}_0)$
are isomorphisms
and, as a consequence, they interchange the only non-zero elements of these two groups
(which are $l_i$ and $l_{i-a}$ if $i<D/2$, or $h^{D-i}$ and $h^{D-i-a}$ if $i>D/2$).
For $i<a$, the projection and the inclusion are isomorphisms between $\CH_i(\bar{X})$ and
$\Z/2\Z=\CH_i(\Z(i))$ and the only non-zero element of the first Chow group is $l_i$.
Finally, for $i>D-a$, the projection and the inclusion are isomorphisms between
$\CH_i(X)$ and $\Z/2\Z$ and the only non-zero element of the Chow group is $h^{D-i}$.

For $i=D/2$ (here we are in the case of even $D$, of course),
the basis of the group $\CH_i(\bar{X})$ is given by the elements
$h^d$ and $l_d$, while
the basis of the group $\CH_{i-a}(\bar{X}_0)$ is given by the elements
$h^{d-a}$ and $l_{d-a}$.
The subgroups $\BCH_d(X)\subset\CH_d(\bar{X})$ and
$\BCH_{d-a}(X_0)\subset\CH_{d-a}(\bar{X}_0)$, however, are $1$-dimensional,
generated by $h^d$ and $h^{d-a}$
(because
$l_{d-a}\not\in\BCH(X_0)$ by Corollary \ref{X is non-essential}).
Since these subgroups are interchanged by the projection and the inclusion, $h^d$
corresponds to $h^{d-a}$.
Now there are only two possibilities for the element of $\CH_d(\bar{X})$ corresponding to
$l_{d-a}\in\CH_{d-a}(\bar{X}_0)$: either this is $l_d$ or this is $l_d+h^d$.
Which one of these two possibilities takes place depends on the construction of the motivic
decomposition $(*)$;
but a given motivic decomposition can be always corrected in such a way
that the first possibility
takes place (one can simply use an automorphism of the variety $X_0$ interchanging
$l_{d-a}$ with $l_{d-a}+h^{d-a}$).
\end{proof}

The ``most important'' summand in the motivic decomposition of $X^r$ is, of course, $X_0^r$.
We introduce a special notation for the projection and the inclusion  related to
this summand:
$\pr^r=\pr_{a\dots a}$ and $\inc^r=\inc_{a\dots a}$.

\begin{cor}
\label{old isotropic}
The (mutually semi-inverse) homomorphisms
$$
\pr^r\!:\CH(\bar{X}^r)\to\CH(\bar{X}_0^r)
\;\;\text{ and }\;\;
\inc^r\!:\CH(\bar{X}_0^r)\to\CH(\bar{X}^r)
$$
(for any $r\geq1$)
map rational cycles to rational cycles;
moreover, the induced homomorphism
$\pr^r\!:\BCH(X^r)\to\BCH(X_0^r)$ is surjective.
\qed
\end{cor}

Now we get an extended version of Corollary \ref{X is non-essential} which reads as
follows:

\begin{cor}
\label{determining witt index}
For an arbitrary quadric $X$ (isotropic or not) and any integer $i$ one has:
$l_i\in\BCH(X)$ if and only if $\IW(X)>i$
(where $\IW(X)=\IW(\phi)$ is the Witt index).
\end{cor}

\begin{proof}
The ``if'' part of the statement is trivial.
Let us prove the ``only if'' part using an induction on $i$.
The case of $i=0$ is served by Lemma \ref{springer}.

Now we assume that $i>0$ and $l_i\in\BCH(X)$.
Since $l_i\cdot h=l_{i-1}$, the element $l_{i-1}$ is rational as
well, and by the induction hypothesis $\IW(X)\geq i$.
If $\IW(X)=i$, then
the image of $l_i\in\BCH(X)$ under the map
$\pr^1\!:\CH(\bar{X})\to\CH(\bar{X}_0)$ equals $l_0$ and is rational
by Corollary \ref{old isotropic}.
Therefore, by Lemma \ref{springer}, the quadric $X_0$ is isotropic,
a contradiction.
\end{proof}

We recall that the {\em splitting pattern}
$\SP(\phi)$
of an anisotropic quadratic form $\phi$
is defined as the set of integers
$$
\SP(\phi)=
\{\dim(\phi_E)_{\an}|\;\text{$E/F$ is a field extension}\}
$$
(here $\phi_E$ stays for the quadratic form over the field $E$ obtained from $\phi$ by
extending the scalars; $(\phi_E)_{\an}$ is the anisotropic part of $\phi_E$).

The splitting pattern can be obtained using the {\em generic
splitting tower} of M. Knebusch (arbitrary filed extensions of $F$ are then replaced by
concrete fields).
To construct this tower, we put $F_1=F(X)$, the
function field of
the projective quadric $X$ given by $\phi$.
Then we put $\phi_1=(\phi_{F_1})_{\an}$ and write $X_1$ for the
projective quadric (over the field $F_1$) given by the quadratic form $\phi_1$.
We proceed by setting $F_2=F_1(X_1)$ and so on until we can
(we stop on $F_\Hight$ such that $\dim(\phi_\Hight)\leq1$).
The tower of fields $F\subset F_1\subset\dots\subset F_\Hight$ obtained
this way is called the generic splitting tower of $\phi$
and (see \cite{Knebusch})
$$
\SP(\phi)=\{\dim(\phi),\dim(\phi_{F_1})_{\an},\dots,\dim(\phi_{F_\Hight})_{\an}\}=
\{\dim(\phi),\dim(\phi_1),\dots,\dim(\phi_\Hight)\}
$$
(the integer $\Hight=\Hight(\phi)$ is the {\em height} of $\phi$;
note that the elements of $\SP(\phi)$ are written down in the descending order).

An equivalent invariant of $\phi$ is called the {\em higher Witt
indices} of $\phi$ and defined as follows.
Let us write the set of integers $\{\IW(\phi_E)|\;\text{$E/F$ a field
extension}\}$, where $\IW(\phi_E)$ is the usual Witt index of
$\phi_E$, in the form
$$
\{0=\IW(\phi)<\iw_1<\iw_1+\iw_2<\dots<\iw_1+\iw_2+\dots+\iw_\Hight\}\;.
$$
The sequence of the positive integers $\iw_1,\dots, \iw_\Hight$ is called the higher
Witt indices of $\phi$.
For every $q\in\{0,1,\dots,\Hight\}$, we also set
$$
\jw_q=\jw_q(\phi)=\iw_0+\iw_1+\dots+\iw_q=\IW(\phi_{F_q})\;.
$$
Clearly, one has
$$
\{0=\jw_0,\jw_1,\dots,\jw_\Hight\}=\{\IW(\phi_E)|\;\text{$E/F$ is a field
extension}\}
$$
(this set of integers is sometimes also called the splitting pattern of
$\phi$ in the literature).

The following
easy observation
is crucial:

\begin{thm}
\label{determining sp}
The splitting pattern as well as the higher Witt indices of
an anisotropic quadratic form $\phi$
(of some given dimension)
are determined by
the group
$$\BCH(X^\ast)=\Oplus_{r\geq1}\BCH(X^r)\;.$$
\end{thm}

\begin{proof}
The pull-back homomorphism
$g_1^*\!:\CH(X^r)\to\CH(X_{F(X)}^{r-1})$ with respect to the
morphism of schemes $g_1\!:X_{F(X)}^{r-1}\to X^r$ given by the generic
point of, say, the first factor of $X^r$, is surjective.
It induces an epimorphism
$\BCH(X^r)\onto\BCH(X_{F(X)}^{r-1})$, which is a restriction of the
epimorphism $\CH(\bar{X}^r)\onto\CH(\bar{X}_{F(X)}^{r-1})$
mapping
each basis element of the shape $h^0\times\beta$ with $\beta\in\CH(\bar{X}^{r-1})$
to $\beta\in\CH(\bar{X}^{r-1}_{F(X)})$ and killing all other basis elements.
Therefore the group $\BCH(X^\ast)$ determines the group
$\BCH(X_{F(X)}^\ast)$.
In particular, the group $\BCH(X_{F(X)})$ is determined, so that we
have reconstructed $\IW(X_{F(X)})=\iw_1(X)$ (see Corollary \ref{determining witt
index}).
Moreover, by Corollary \ref{old isotropic}, the group $\BCH(X_{F(X)}^\ast)$
determines the group $\BCH(X_1^\ast)$
(via the surjection $\pr^\ast$; $X_1$ staying for the anisotropic part of
$X_{F(X)}$), and we can proceed by induction.
\end{proof}

\begin{rem}
The proof of Theorem \ref{determining sp} makes it clear that the
statement of this theorem can be made more precise in the
following way.
If for some $q\in \{1,2,\dots,\Hight\}$ the Witt indices $\iw_1,\dots,\iw_{q-1}$ are
already reconstructed, then one determines
$\iw_q=\jw_q-\iw_{q-1}-\dots-\iw_1$ by
the formula
$$
\jw_q=\max\{j|\;
\text{the product $h^0\times h^{\jw_1}\times\dots\times h^{\jw_{q-1}}\times l_{j-1}$
is contained in a rational cycle}\}\;.
$$
\end{rem}

\begin{rem}
Concluding this section, we would like to underline that the role
of the algebraic closure $\bar{F}$ in the definition of the group
$\BCH(X^\ast)$ is secondary:
the group $\CH(\bar{X}^\ast)$ (used in the definition of
$\BCH(X^\ast)$) has to be interpreted as the direct limit
$\lim\CH(X^\ast_E)$ taken over {\em all} field extensions
$E/F$.
The homomorphism $\CH(\bar{X}^\ast)\to\lim\CH(X^\ast_E)$ is
an isomorphism.
More generally, the homomorphism
$\CH(X^\ast_E)\to\lim\CH(X^\ast_E)$ for a given $E/F$ is an
isomorphism if and only if the quadratic form $\phi_E$ is
completely split
(in particular, for any $E/F$ with completely split $\phi_E$,
there is a canonical isomorphism
$\CH(X_E^\ast)=\CH(\bar{X}^\ast)$, coinciding with the
composition
$\res^{-1}_{\bar{E}/\bar{F}}\compose\res_{\bar{E}/E}$,
where $\bar{E}$ is a field containing $E$ and $\bar{F}$).
\end{rem}

\section
{Strategy of proof}
\label{Strategy of proof}

As shown in Theorem \ref{determining sp},
the group $\BCH(X^\ast)$ determines
the splitting pattern of
the quadratic form $\phi$.
In its turn, the splitting pattern of $\phi$ determines the powers
of the fundamental ideal of the Witt ring containing the class of
$\phi$.
At least it is easy to prove

\begin{lemma}
\label{sp-I^n}
Let $\phi$ be an even-dimensional anisotropic quadratic form and let $n\geq1$ be an
integer. We write $p$ for the least positive integer of
$\SP(\phi)$ (note that $p$ is a power of $2$, \cite[thm. 5.4(i)]{Scharlau}).
If $[\phi]\in I^n$, then $p\geq 2^n$.
\end{lemma}

\begin{proof}
We assume that $[\phi]\in I^n(F)$.
Let $E/F$ be a field extension such that $\dim(\phi_E)_{\an}=p$.
Since $0\ne[(\phi_E)_{\an}]\in I^n(E)$, we get that $p\geq 2^n$
by the Arason-Pfister theorem.
\end{proof}

\begin{rem}
It is not needed in this paper but nevertheless good to know that
the converse statement to Lemma \ref{sp-I^n} is also true.
This is a hard result, however.
It is proved in \cite[thm. 4.3]{OVV}.
\end{rem}

Now we are able to describe the strategy of our proof of
Conjecture \ref{conj}.
Let us consider a power $I^n$ of the fundamental ideal.
Let $\phi$ be an anisotropic quadratic form with $[\phi]\in I^n$,
having some dimension prohibited by Conjecture \ref{conj}.
The group $\BCH(X^\ast)$, where $X$ is the projective quadric
given by $\phi$, should satisfy some restrictions listed bellow.
This group is a subgroup of $\CH(\bar{X}^\ast)$, the later one depends
only on the dimension of $\phi$.
So, we prove Conjecture \ref{conj}, if we check that every subgroup
of $\CH(\bar{X}^\ast)$, satisfying the list of restrictions, can
not be $\BCH(X^\ast)$ for the form $\phi$
(see Lemma \ref{sp-I^n} with Theorem \ref{determining sp}).
This is the way we prove Conjecture \ref{conj}.

And here is the list of restrictions on $\BCH(X^\ast)$ considered as a subset
of $\CH(\bar{X}^\ast)$ (a big part of this list is of course valid for an
arbitrary smooth projective variety in place of the quadric $X$):

\begin{prop}
\label{list}
Assuming that the quadric $X$ is anisotropic,\footnote{Anisotropy is important only
for (\ref{springersatz}), (\ref{binary size}), (\ref{9}), and (\ref{10}).} we have:
\begin{enumerate}
\item
\label{addition and multiplication}
$\BCH(X^\ast)$ is closed under
addition and multiplication;
\item
$\BCH(X^\ast)$ is closed under
passing to the homogeneous components (with respect to the grading
of the Chow group and to the $\ast$-grading);
\item
\label{h^0 and h^1}
$\BCH(X^\ast)$ contains $h^0=[\bar{X}]$ and $h^1=h$
(and therefore contains $h^i$ for all $i\geq0$);
\item
\label{springersatz}
{\em [{\bf Springer-Satz}]}\\
$\BCH(X^\ast)$
does not contain $l_0$;
\item
for every $r\geq1$,
$\BCH(X^r)$ is closed under the automorphisms of $\CH(\bar{X}^r)$
given by the permutations of factors of $X^r$;
\item
\label{projections and diagonals}
for every $r\geq1$,
$\BCH(X^\ast)$ is closed under push-forwards and pull-backs
with respect to all $r$ projections $X^r\to X^{r-1}$ and to all
$r$ diagonals $X^r\to X^{r+1}$
(taking in account the previous restriction, it is enough to speak only about
the first projection
$$x_1\times x_2\times\dots\times x_r\mapsto x_2\times\dots\times x_r
$$
and the first diagonal
$$x_1\times x_2\times\dots\times x_r\mapsto x_1\times x_1\times x_2\dots\times x_r
$$
here);
\item
$\BCH(X^\ast)$ is closed under the total Steenrod operation
$$
S\!:\CH(\bar{X}^\ast)\to\CH(\bar{X}^\ast)\;;
$$
\item
\label{binary size}
{\em[A. Vishik, {\bf ``Size of binary correspondences''}]}\\
if
$
\BCH(X^2)\ni h^0\times l_i+l_i\times h^0
$
for some integer $i\geq0$, then
the integer $\dim(X)-i+1$ is a power of $2$;
\item
\label{9}
{\em [{\bf ``Inductive restriction''}]}\\
the image of $\BCH(X^{\ast+1})$ under the composition
$$
\begin{CD}
\CH(\bar{X}^{\ast+1})@>{g_1^\ast}>>\CH(\bar{X}_{F(X)}^\ast)@>{\pr^\ast}>>
\CH(\bar{X}_1^\ast)
\end{CD}
$$
($g_1^\ast$ is introduced in the proof of Theorem
\ref{determining sp},
$\pr^\ast$ in Corollary \ref{old isotropic})
should coincide with $\BCH(X_1^\ast)$ and therefore should
satisfy all restrictions listed in this proposition (including the current
one);
\item
\label{10}
{\em [{\bf ``Supplement to inductive restriction''}]}\\
for any integer $r\geq2$, any integer $s\in[1,\;r)$, and any
projection $\pr_{i_1\dots i_r}$ of Proposition \ref{isotropic},
the image of $\BCH(X^r)$ under
$\pr_{i_1\dots i_r}\!:\CH(\bar{X}^r)\to\CH(\bar{X}_1^{s(i_1,\dots,i_r)})$
is inside of
$\BCH(X_1^{s(i_1,\dots,i_r)})$ (reconstructed by {\em (\ref{9})}).
\end{enumerate}
\end{prop}

\begin{proof}
Only the property (\ref{binary size}) needs a proof.
We note that this property does not seem to be a consequence of
the others.
It is proved in \cite[thm. 5.1]{KMD} by some computation in the {\em integral}
Chow group of $X^\ast$, not in the modulo $2$ Chow group (although involving
modulo $2$ Steenrod operations)
(the original proof is in \cite[thm. 6.1]{Izhboldin-Vishik}; it makes use of
higher motivic cohomology).

More precisely, the case of $i=0$ is proved in \cite[thm.
5.1]{KMD}.
In order to reduce the general case to the case of $i=0$,
we take an arbitrary subquadric $Y\subset X$ of
codimension $i$ and pullback the cycle $h^0\times l_i+l_i\times h^0$
with respect to the imbedding $Y^2\hookrightarrow X^2$.
The result is $h^0\times l_0+l_0\times h^0\in\BCH(Y^2)$.
Therefore, $\dim(Y)+1$ is a power of $2$ by \cite[thm. 5.1]{KMD}.
Since $\dim(Y)=\dim(X)-i$,
it follows that the integer $\dim(X)-i+1$ is the same power of
$2$.
\end{proof}

\begin{rem}
Obviously, one can write down some additional restrictions on
$\BCH(X^\ast)$.
However, all restrictions I know are consequences of the
restrictions of Proposition \ref{list}.
For instance, $\BCH(X^\ast)$ should be stable with respect to
the external products; but this is a consequence of the stability
with respect to the internal products (\ref{addition and multiplication})
and the pull-backs with respect to projections (\ref{projections and
diagonals}).
Another example:
the image of the total Chern class
$c\!:K_0(\bar{X}^\ast)\to\CH(\bar{X}^\ast)$
restricted to $K_0(X^\ast)$
(note that $K_0(X^\ast)$ is computed for quadrics \cite{Swan} and, more
generally, for all projective homogeneous varieties \cite{Panin})
should be inside of $\BCH(X^\ast)$;
but it is already guarantied by the fact that $\BCH(X^\ast)$ is
closed under addition and multiplication (\ref{addition and
multiplication}) and contains $h^0$ and $h^1$
(\ref{h^0 and h^1}).\footnote{However
the property with the Chern class can be a good replacement
for (\ref{h^0 and h^1}) when
transfering this theory to other algebraic varieties in place of the quadric $X$.}
One more example:
$\BCH(X^2)$ should be closed under  the composition of
correspondences (see \cite[\S16]{Fulton} for the definition of composition of
correspondences);
but the operation of composition of correspondences is produced by
pull-backs and push-forwards with respect to projections together
with the operation of multiplication of cycles.
\end{rem}

\begin{rem}
\label{computability}
Let us remark that all operations involved in the list of
restrictions of Proposition \ref{list} are easy to compute in terms of the basis
elements.
The multiplication in $\CH(\bar{X}^\ast)$ was
described in the previous {\S};
a formula for the total Steenrod operation was given already as
well.
Also the operations used in the inductive restriction are computed
(see Corollary \ref{old isotropic} and the proof of Theorem
\ref{determining sp}).
As to the pull-backs and push-forwards with respect to the first
projection $\pr\!:X^{r+1}\to X^r$ and to the first diagonal
$\delta\!:X^r\to X^{r+1}$, they are computed for basis elements
$\beta_0,\beta_1,\dots,\beta_r\in\CH(\bar{X})$ as follows:
$$
\begin{array}{lll}
\pr^*(\beta_1\times\dots\times\beta_r)&=
&h^0\times\beta_1\times\dots\times\beta_r\;;\\[1em]
\pr_*(\beta_0\times\beta_1\times\dots\times\beta_r)&=&
\begin{cases}
\beta_1\times\dots\times\beta_r\;,\;&\text{if $\beta_0=l_0$,}\\
0\;,&\text{otherwise;}
\end{cases}\\[2em]
\delta^*(\beta_0\times\beta_1\times\dots\times\beta_r)&=&
(\beta_0\cdot\beta_1)\times\beta_2\times\dots\times\beta_r\;;\\[1em]
\delta_*(\beta_1\times\dots\times\beta_r)&=&
\big((\beta_1\times h^0)\cdot\Delta\big)
\times\beta_2\times\dots\times\beta_r
=\\
&&\big((h^0\times\beta_1)\cdot\Delta\big)
\times\beta_2\times\dots\times\beta_r
\;;
\end{array}
$$
where $\Delta\in\BCH(X^2)$ is the class of the diagonal
computed in Corollary \ref{diagonal}.
\end{rem}

\begin{rem}
One can obviously simplify a little bit the
list of restrictions of Proposition \ref{list}.
For instance, instead of stability under the push-forwards with
respect to the diagonals, it suffices to require that the
cycle $\Sum_{i=0}^d (h^i\times l_i+l_i\times h^i)$,
related to the diagonal,
lies in $\BCH(X^2)$ (see Remark \ref{computability} and Corollary \ref{diagonal}).
Also the inductive restriction is not so restrictive as it may
seem to: the group $\BCH(X_1^\ast)$ automatically satisfies
most of the required restrictions.
\end{rem}

\begin{rem}
Looking at the list of restrictions, it is easy to see that every
group $\BCH(X^r)$ determines $\BCH(X^{<r})$.
Moreover, one can show that $\BCH(X^d)$ determines the whole
group $\BCH(X^\ast)$.~\footnote{It would be interesting
to rewrite all restrictions of
Proposition \ref{list} in terms of the group $\BCH(X^d)$.
}
Since $\BCH(X^d)$ is a subgroup of the finite group
$\CH(\bar{X}^d)$, it follows, in particular, that the invariant
$\BCH(X^d)$ (of the quadratic forms $\phi$ of a given dimension)
has only a finite number of different values
(this way one also sees that Conjecture \ref{conj} can
be checked for any
concrete dimension by computer).
\end{rem}

We will often use the composition of correspondences,
even for the cycles on bigger than $2$ powers of $X$: this is
a convenient way to handle the things.
Namely, for $\alpha\in\BCH(X^r)$ and
$\alpha'\in\BCH(X^{r'})$, we may consider $\alpha$ as a
correspondence, say, from $X^{r-1}$ to $X$, and we may consider
$\alpha'$ as a correspondence from $X$ to $X^{r'-1}$; then the
composition $\alpha'\compose\alpha$ is a cycle in
$\BCH(X^{r+r'-2})$,
and here is the formula for composing the basis elements
(we put here this obvious formula because it will be used many times in our computations):

\begin{lemma}
\label{comp-corr}
The composition $\beta'\compose\beta\in\CH(\bar{X}^{r+r'-2})$ of two basis elements
$$
\beta=\beta_1\times\dots\times\beta_r\in\CH(\bar{X}^r)\;\;
\text{ and }\;\;
\beta'=\beta'_1\times\dots\times\beta'_{r'}\in\CH(\bar{X}^{r'})
$$
is equal to
$\beta_1\times\dots\times\beta_{r-1}\times\beta'_2\times\dots\times\beta'_{r'}$,
if $\beta_r\cdot\beta'_1=l_0$;
otherwise the composition $\beta'\compose\beta$
is $0$.
\qed
\end{lemma}

\begin{cor}
\label{diagonal}
For the diagonal class $\Delta\in\BCH(X^2)$, one has:
$$
\Delta=\Sum_{i=0}^d(h^i\times l_i+l_i\times
h^i)+(D+1)\cdot(d+1)\cdot (h^d\times h^d)\;.
$$
In particular, the sum $\Sum_{i=0}^d(h^i\times l_i+l_i\times
h^i)$ is always rational.
\end{cor}

\begin{proof}
Using Lemma \ref{comp-corr}, it is straight-forward to verify
that the cycle given by the above
formula acts (by composition) trivially on any basis cycle of
$\CH^2(\bar{X}^2)$.
\end{proof}

\section
{Cycles on $X^2$}

We are using the notation introduced in
{\S} \ref{Notation and preliminary observations}.
In particular,
$X$ is the projective quadric of an even dimension
$\D=2d$ or of an odd dimension $\D=2d+1$ given by a quadratic form $\phi$ over the field $F$.
In this {\S} we assume that $X$ is anisotropic.
Let $\iw_1,\dots \iw_\Hight$ be the higher Witt indices of $\phi$ (with $\Hight$ being the height of
$\phi$).
We write $S$ for the set $\{1,2,\dots,\Hight\}$ and we set
$\jw_q=\iw_1+\iw_2+\dots+\iw_q$ for every $q\in S$.

An ``important'' as well as the ``first interesting'' part
of the group $\BCH(X^\ast)$ is the group
$\BCH(X^2)$ and especially $\BCH_D(X^2)=\BCH^D(X^2)$
(note that, due to the Rost nilpotence theorem
(\cite{Rost}, see also \cite{Brosnan-nilpotence}),
the latter group detects all motivic decompositions of $X$).
The group $\BCH_D(X^2)$ was studied intensively by A. Vishik
(see, e.g., \cite{Vishik-Lens}).
In the next {\S} we reproduce most of his results concerning this
group.
%
Actually we give an extended version of these results
describing the structure of a bigger group, namely of the group $\BCH^{\leq D}(X^2)$.

Originally, Vishik's results are formulated in terms of motivic
decompositions of $X$;
by this reason, their proofs use the Rost nilpotence theorem for
quadrics,
which is not used in the present text at all.
Here we
simplify the formulation;
we also give a different
complete self-contained
proof and show that
all these results are consequences of the restrictions on
$\BCH(X^\ast)$ listed in Proposition \ref{list}.
More precisely, we start with some general results concerning
an arbitrary anisotropic quadric $X$;
the proofs of these results
use neither the restriction provided by the Steenrod
operation, nor the ``size of binary correspondences'' restriction;
a summary of these results is given in Theorem \ref{kruto!}.
Proposition \ref{known}, appearing in the very end of the section, contains
the result on the so called small quadrics (Definition
\ref{def-small}), needed in the proof of Conjecture \ref{conj};
its proof uses the ``size of binary correspondences'' restriction
(the Steenrod operation does still not show up).

We start by the following easy observation:

\begin{lemma}
\label{herovyj element}
For an even $D=\dim(X)$,
the basis element $l_d\times l_d\in\BCH_D(X^2)$ does
not appear in the decomposition of a rational cycle.
\end{lemma}

\begin{proof}
We assume the contrary.
Let $\alpha$ be a cycle in $\BCH_D(X^2)$ containing $l_d\times l_d$.
Then the push-forward with respect to the projection  onto the
second factor
$\pr\!:X^2\to X$ of the cycle $\alpha\cdot(h^d\times h^0)$ is rational
and equals $l_d$ or $h^d+l_d$
(because $\pr_*\big(\beta\cdot(h^d\times h^0)\big)$ is $l_d$ for
$\beta=(l_d\times l_d)$, $h^d$ for $\beta=l_d\times h^d$, and $0$ for every other
basis cycle $\beta\in\CH_D(\bar{X}^2)$).
Therefore the cycle $l_d$ is rational, showing that $X$ is hyperbolic
(Corollary \ref{determining witt index}), a contradiction.
\end{proof}

\begin{lemma}
\label{intersection}
If $\alpha_1,\alpha_2\in\BCH^{\leq D}(X^2)$,
then  the cycle
$\alpha_1\cap\alpha_2$
is rational (where the notation $\alpha_1\cap\alpha_2$ means the sum of the basis cycles
contained simultaneously in $\alpha_1$ and in $\alpha_2$).
\end{lemma}

\begin{proof}
Clearly, we may assume that $\alpha_1$ and $\alpha_2$ are
homogeneous of the same dimension $D+i$ and do not contain any non-essential basis element.
Then the intersection
(modulo the non-essential elements)
is computed as
$
\alpha_1\cap\alpha_2=\alpha_2\compose\Big(\alpha_1\cdot(h^0\times h^i)\Big)
$
(see Lemma \ref{comp-corr}).
\end{proof}

\begin{dfn}
We write $\BCHE(\bar{X}^\ast)$ for the subgroup of
$\CH(\bar{X}^\ast)$ generated by the essential basis elements.
We set
$
\BCHE(X^\ast)=\BCHE(\bar{X}^\ast)\cap\BCH(X^\ast)\;.
$
Note that the group $\BCHE(X^\ast)$ is a subgroup of $\BCH(X^\ast)$
isomorphic to the
quotient of $\BCH(X^\ast)$ by the subgroup of the non-essential
elements.
\end{dfn}

\begin{dfn}
A non-zero cycle of $\BCH^{\leq D}(X^2)$ is called {\em minimal}, if it does not contain
any proper rational subcycle.
Note that a minimal cycle does not contain non-essential basis
elements, i.e., lies in $\BCHE(X^2)$.
Also note that a minimal cycle is always homogeneous.
\end{dfn}

A very first structure result on $\BCH^{\leq D}(X^2)$ reads as
follows:

\begin{prop}
\label{minimal form basis}
The minimal cycles form a basis of the group $\BCHE^{\leq D}(X^2)$.
Two different minimal cycles do not intersect each other
(here we speak about the notion of intersection of cycles introduced in Lemma
\ref{intersection}).
The sum of the minimal cycles of dimension $D$
is equal to the
sum
$
\Sum_{i=0}^d h^i\times l_i+l_i\times h^i
$
of all $D$-dimensional essential basis elements (excluding $l_d\times l_d$
in the case of even $D$).
\end{prop}

\begin{proof}
The first two statements of Proposition \ref{minimal form basis}
follow from Lemma \ref{intersection}.
The third statement follows from previous ones together with the
rationality of the diagonal cycle
(see Corollary \ref{diagonal}).
\end{proof}

\begin{dfn}
Let $\alpha$ be a homogeneous cycle in $\CH^{\leq D}(\bar{X}^2)$.
For every $i$ with $0\leq i\leq \dim(\alpha)-D$, the products
$\alpha\cdot(h^0\times h^i)$,
$\alpha\cdot(h^1\times h^{i-1})$, and so on up to
$\alpha\cdot(h^i\times h^0)$ will be called the ($i$-th order)
{\em derivatives} of $\alpha$.
Note that all derivatives are still in $\CH^{\leq D}(\bar{X}^2)$
and that all derivatives of a rational cycle are also rational.
\end{dfn}

\begin{lemma}
\label{trivial on derivatives}
\begin{enumerate}
\item
Each derivative of any essential basis element $\beta\in\BCHE^{\leq D}(\bar{X}^2)$
is an essential basis element.
\item
For any $r\geq0$ and any essential basis cycles
$\beta_1,\beta_2\in\BCHE_{D+r}(\bar{X}^2)$, two derivatives
$
\beta_1\cdot(h^{i_1}\times h^{j_1})\;\;\text{ and }\;\;
\beta_2\cdot(h^{i_2}\times h^{j_2})
$
of $\beta_1$ and $\beta_2$ coincide only if
$\beta_1=\beta_2$, $i_1=i_2$, and $j_1=j_2$.
\end{enumerate}
\end{lemma}

\begin{proof}
(1)
If $\beta\in\BCHE_{D+r}(\bar{X}^2)$ for some $r\geq0$, then,
up to transposition, $\beta=h^i\times l_{i+r}$ with
$i\in[0,\;d-r]$.
An arbitrary derivative of $\beta$ is equal to
$
\beta\cdot(h^{j_1}\times h^{j_2})=h^{i+j_1}\times l_{i+r-j_2}
$
with some $j_1,j_2\geq0$ such that
$j_1+j_2\leq r$.
We have $0\leq i+j_1\leq d$ and therefore $h^{i+j_1}$ is a basis
element.
We also have $d\geq i+r-j_2\geq0$ and therefore $l_{i+r-j_2}$ is a
basis element too.

Statement (2) is
trivial.
\end{proof}

\begin{rem}
\label{small pyramid}
For the sake of visualization, it is good to think of the basis cycles
of $\CH^{\leq D}(\bar{X}^2)$ (with $l_{D/2}\times l_{D/2}$ excluded by the reason of
Lemma \ref{herovyj element}) as of points of the ``pyramid''
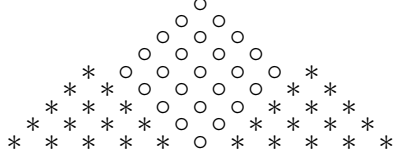
$$
\renewcommand{\arraycolsep}{0.1ex}
\renewcommand{\arraystretch}{0.25}
\begin{array}{ccccccccccccccccccccc}
&&&&&&&&&&\NBC&&&&&&&&&&\\
&&&&&&&&&\NBC&&\NBC&&&&&&&&&\\
&&&&&&&&\NBC&&\NBC&&\NBC&&&&&&&&\\
&&&&&&&\NBC&&\NBC&&\NBC&&\NBC&&&&&&&\\
&&&&\BC&&\NBC&&\NBC&&\NBC&&\NBC&&\NBC&&\BC  &&&&\\
&&&\BC&&\BC&&\NBC&&\NBC&&\NBC&&\NBC&&\BC&&\BC  &&&\\
&&\BC&&\BC&&\BC&&\NBC&&\NBC&&\NBC&&\BC&&\BC&&\BC  &&\\
&\BC&&\BC&&\BC&&\BC&&\NBC&&\NBC&&\BC&&\BC&&\BC&&\BC &\\
\BC&&\BC&&\BC&&\BC&&\BC&&\NBC&&\BC&&\BC&&\BC&&\BC&&\BC
\end{array}
$$
($D=8$ on the picture;
for an odd $D$ the pyramid has no ``step'', see, e.g., the picture of Remark
\ref{big piramid}), where the $\NBC$-s stay for the non-essential
basis elements while the $\BC$-s stay for the essential ones;
the point on the top is $h^0\times h^0$;
for every $i\in\{0,1,\dots,D\}$, the $i$-th row of the pyramid represents the basis of
$\CH^i(\bar{X}^2)$ (in the case of even $D$, the $D$-th row is the basis without
$l_d\times l_d$) ordered by the codimension of the first factors
(starting with $h^0\times?$).
For any $\alpha\in\BCHE^{\leq D}(\bar{X}^2)$, we can put a mark on
the points representing basis elements contained in the
decomposition of $\alpha$; the set of marked points is the {\em
diagram} of $\alpha$.
If $\alpha$ is homogeneous, the marked points lie in the same
row.
Now it is easy to interpret the derivatives of $\alpha$:
the diagram of an $i$-th order derivative is a projection of the marked
points of the diagram of $\alpha$ to the $i$-th row bellow along
some direction.
The diagram of every derivative of $\alpha$ has the same number
of marked points as the diagram of $\alpha$ (Lemma
\ref{trivial on derivatives}).
The diagrams of two different derivatives of the same order are shifts
(to the right or to the left) of each other.
\end{rem}

\begin{lemma}
\label{4.9}
The following conditions on a homogeneous cycle $\alpha\in\BCH^{\leq D}(X^2)$
are equivalent:
\begin{enumerate}
\item
$\alpha$ is minimal;
\item
all derivatives of $\alpha$ are minimal;
\item
at least one derivative of $\alpha$ is minimal.
\end{enumerate}
\end{lemma}

\begin{proof}
Derivatives of a proper subcycle of $\alpha$ are proper subcycles
of the derivatives of $\alpha$; therefore, $(3)\Rightarrow(1)$.

In order to show that $(1)\Rightarrow(2)$, it suffices to show
that
two first order derivatives
$\alpha\cdot(h^0\times h^1)$ and $\alpha\cdot(h^1\times h^0)$
of a minimal cycle $\alpha$ are minimal.
In the contrary case,
replacing eventually $\alpha$ by its transposition, we come to
the situation where the derivative $\alpha\cdot(h^0\times h^1)$
of a minimal $\alpha$ is not minimal.
It follows that the cycle $\alpha\cdot(h^0\times h^i)$,
where $i=\dim(\alpha)-D$,
is not minimal too;
let $\alpha'$ be its proper subcycle.
Taking the composition $\alpha\compose\alpha'$ and removing the non-essential
summands,
we get a proper subcycle
of $\alpha$ (see Lemma \ref{comp-corr}).
\end{proof}

\begin{cor}
The derivatives of a minimal cycle are disjoint.
\end{cor}

\begin{proof}
The derivatives of a minimal cycle are minimal (Lemma \ref{4.9})
and pairwise different (Lemma \ref{trivial on derivatives}).
Two different minimal cycles are disjoint by Lemma \ref{intersection}
(see also Proposition \ref{minimal form basis}).
\end{proof}

\begin{lemma}
\label{dyrki v verhnem cykle}
Let $\alpha$ be an element of $\BCH_{D+k-1}(X^2)$ with some
$k\geq1$.
For any $q\in\{1,\dots,\Hight\}$ and for any non-negative $i$ with $\iw_q-k< i<\iw_q$,
the cycle
$\alpha$ contains neither the product $h^{\jw_{q-1}+i}\times
l_{\jw_{q-1}+i+k-1}$ nor the transposition  of this product.
\end{lemma}

\begin{proof}
Let us assume the contrary:
for some $k\geq1$,
some $q\in\{1,\dots,\Hight\}$, and some $i$ with $\iw_q-k< i<\iw_q$,
there exists a rational cycle $\alpha$
containing the product $h^{\jw_{q-1}+i}\times
l_{\jw_{q-1}+i+k-1}$ or the transposition of this product.
If $\alpha$ contains the transposition of the product, we replace
$\alpha$ by the transposition of $\alpha$.
Passing to the ($q-1$)-th field of the generic splitting tower and
using the projection of Corollary \ref{old isotropic},
we come to the situation where $q=1$ and
$\alpha$ contains the product $h^i\times l_{i+k-1}$
such that $\iw_1-k<i<\iw_1$.
The projection
$\pr_{D-i,\iw_1}(\alpha_{F(X)})$ is a rational cycle on $\bar{X_1}$ containing
$l_{i+k-1-\iw_1}$ (note that $i+k-1-\iw_1\geq0$).
We get a contradiction with Corollary \ref{X is non-essential}.
\end{proof}

\begin{rem}
\label{big piramid}
In order to ``see'' the statement of Lemma \ref{dyrki v verhnem
cykle}, it is helpful to mark by $\AC$ the essential basis elements
which are not ``forbidden'' by this lemma
(we are speaking about the pyramid of basis cycles drawn in Remark
\ref{small pyramid}).
We will get isosceles triangles based on the lower row of the pyramid.
For example, if $X$ is a $29$-dimensional quadric with the higher
Witt indices $4,3,5,2$ (such a quadric $X$ does not
exist in reality, but is convenient for the illustration),
then the picture looks as follows:

$$
\renewcommand{\arraycolsep}{0.005ex}
\renewcommand{\arraystretch}{0.0125}
\begin{array}{ccccccccccccccccccccccccccccccccccccccccccccccccccccccc}
&&&&&&&&&&&&&&&&&&&&&&&&&&&\NC&&&&&&&&&&&&&&&&&&&&&&&&&&&\\
&&&&&&&&&&&&&&&&&&&&&&&&&&\NC&&\NC&&&&&&&&&&&&&&&&&&&&&&&&&&\\
&&&&&&&&&&&&&&&&&&&&&&&&&\NC&&\NC&&\NC&&&&&&&&&&&&&&&&&&&&&&&&&\\
&&&&&&&&&&&&&&&&&&&&&&&&\NC&&\NC&&\NC&&\NC&&&&&&&&&&&&&&&&&&&&&&&&\\
&&&&&&&&&&&&&&&&&&&&&&&\NC&&\NC&&\NC&&\NC&&\NC&&&&&&&&&&&&&&&&&&&&&&&\\
&&&&&&&&&&&&&&&&&&&&&&\NC&&\NC&&\NC&&\NC&&\NC&&\NC&&&&&&&&&&&&&&&&&&&&&&\\
&&&&&&&&&&&&&&&&&&&&&\NC&&\NC&&\NC&&\NC&&\NC&&\NC&&\NC&&&&&&&&&&&&&&&&&&&&&\\
&&&&&&&&&&&&&&&&&&&&\NC&&\NC&&\NC&&\NC&&\NC&&\NC&&\NC&&\NC&&&&&&&&&&&&&&&&&&&&\\
&&&&&&&&&&&&&&&&&&&\NC&&\NC&&\NC&&\NC&&\NC&&\NC&&\NC&&\NC&&\NC&&&&&&&&&&&&&&&&&&&\\
&&&&&&&&&&&&&&&&&&\NC&&\NC&&\NC&&\NC&&\NC&&\NC&&\NC&&\NC&&\NC&&\NC&&&&&&&&&&&&&&&&&&\\
&&&&&&&&&&&&&&&&&\NC&&\NC&&\NC&&\NC&&\NC&&\NC&&\NC&&\NC&&\NC&&\NC&&\NC&&&&&&&&&&&&&&&&&\\
&&&&&&&&&&&&&&&&\NC&&\NC&&\NC&&\NC&&\NC&&\NC&&\NC&&\NC&&\NC&&\NC&&\NC&&\NC&&&&&&&&&&&&&&&&\\
&&&&&&&&&&&&&&&\NC&&\NC&&\NC&&\NC&&\NC&&\NC&&\NC&&\NC&&\NC&&\NC&&\NC&&\NC&&\NC&&&&&&&&&&&&&&&\\
&&&&&&&&&&&&&&\NC&&\NC&&\NC&&\NC&&\NC&&\NC&&\NC&&\NC&&\NC&&\NC&&\NC&&\NC&&\NC&&\NC&&&&&&&&&&&&&&\\
&&&&&&&&&&&&&\BC&&\NC&&\NC&&\NC&&\NC&&\NC&&\NC&&\NC&&\NC&&\NC&&\NC&&\NC&&\NC&&\NC&&\BC&&&&&&&&&&&&&\\
&&&&&&&&&&&&\BC&&\BC&&\NC&&\NC&&\NC&&\NC&&\NC&&\NC&&\NC&&\NC&&\NC&&\NC&&\NC&&\NC&&\BC&&\BC&&&&&&&&&&&&\\
&&&&&&&&&&&\BC&&\BC&&\BC&&\NC&&\NC&&\NC&&\NC&&\NC&&\NC&&\NC&&\NC&&\NC&&\NC&&\NC&&\BC&&\BC&&\BC&&&&&&&&&&&\\
&&&&&&&&&&\BC&&\BC&&\BC&&\BC&&\NC&&\NC&&\NC&&\NC&&\NC&&\NC&&\NC&&\NC&&\NC&&\NC&&\BC&&\BC&&\BC&&\BC&&&&&&&&&&\\
&&&&&&&&&\BC&&\BC&&\BC&&\BC&&\BC&&\NC&&\NC&&\NC&&\NC&&\NC&&\NC&&\NC&&\NC&&\NC&&\BC&&\BC&&\BC&&\BC&&\BC&&&&&&&&&\\
&&&&&&&&\BC&&\BC&&\BC&&\BC&&\BC&&\BC&&\NC&&\NC&&\NC&&\NC&&\NC&&\NC&&\NC&&\NC&&\BC&&\BC&&\BC&&\BC&&\BC&&\BC&&&&&&&&\\
&&&&&&&\BC&&\BC&&\BC&&\BC&&\BC&&\BC&&\BC&&\NC&&\NC&&\NC&&\NC&&\NC&&\NC&&\NC&&\BC&&\BC&&\BC&&\BC&&\BC&&\BC&&\BC&&&&&&&\\
&&&&&&\BC&&\BC&&\BC&&\BC&&\BC&&\BC&&\BC&&\BC&&\NC&&\NC&&\NC&&\NC&&\NC&&\NC&&\BC&&\BC&&\BC&&\BC&&\BC&&\BC&&\BC&&\BC&&&&&&\\
&&&&&\BC&&\BC&&\BC&&\BC&&\BC&&\BC&&\BC&&\BC&&\BC&&\NC&&\NC&&\NC&&\NC&&\NC&&\BC&&\BC&&\BC&&\BC&&\BC&&\BC&&\BC&&\BC&&\BC&&&&&\\
&&&&\BC&&\BC&&\BC&&\BC&&\BC&&\BC&&\BC&&\AC&&\BC&&\BC&&\NC&&\NC&&\NC&&\NC&&\BC&&\BC&&\AC&&\BC&&\BC&&\BC&&\BC&&\BC&&\BC&&\BC&&&&\\
&&&\AC&&\BC&&\BC&&\BC&&\BC&&\BC&&\BC&&\AC&&\AC&&\BC&&\BC&&\NC&&\NC&&\NC&&\BC&&\BC&&\AC&&\AC&&\BC&&\BC&&\BC&&\BC&&\BC&&\BC&&\AC&&&\\
&&\AC&&\AC&&\BC&&\BC&&\AC&&\BC&&\BC&&\AC&&\AC&&\AC&&\BC&&\BC&&\NC&&\NC&&\BC&&\BC&&\AC&&\AC&&\AC&&\BC&&\BC&&\AC&&\BC&&\BC&&\AC&&\AC&&\\
&\AC&&\AC&&\AC&&\BC&&\AC&&\AC&&\BC&&\AC&&\AC&&\AC&&\AC&&\BC&&\AC&&\NC&&\AC&&\BC&&\AC&&\AC&&\AC&&\AC&&\BC&&\AC&&\AC&&\BC&&\AC&&\AC&&\AC&\\
\AC&&\AC&&\AC&&\AC&&\AC&&\AC&&\AC&&\AC&&\AC&&\AC&&\AC&&\AC&&\AC&&\AC&&\AC&&\AC&&\AC&&\AC&&\AC&&\AC&&\AC&&\AC&&\AC&&\AC&&\AC&&\AC&&\AC&&\AC
\end{array}
$$
\end{rem}

\begin{dfn}
\label{def-shell triangles}
The triangles of Remark \ref{big piramid} will be called the
{\em shell triangles} (their bases are shells in the sense of A.
Vishik).
The shell triangles in the left half of the pyramid are counted
from the left starting by $1$.
The shells triangles in the right half of the pyramid are counted
from the right starting by $1$ as well (so that the symmetric triangles have the
same number; for any $q\in S$, the bases of the $q$-th triangles have (each) $\iw_q$
points).
The rows of the shell triangles are counted from bellow starting by
$0$.
The points of rows of the shell triangles (of the left ones as well
as of the right one) are counted from the left starting by $1$.
\end{dfn}

\begin{lemma}
\label{chetnost' chisla tochek}
For every rational cycle $\alpha\in\BCH^{\leq D}(X^2)$, the number
of the essential basis cycles contained in $\alpha$ is even
(that is, the number of the marked points in the diagram of any
$\alpha\in\BCHE^{\leq D}(X^2)$ is even).
\end{lemma}

\begin{proof}
We may assume that $\alpha$ is homogeneous, say,
$\alpha\in\BCH_{D+k}(X^2)$, $k\geq0$.
Let $n$ be the number of the essential basis cycles contained in
$\alpha$.
The pull-back $\delta^*(\alpha)$ of $\alpha$ with respect to the diagonal
$\delta\!:X\to X^2$ produces $n\cdot l_k\in\BCH(X)$.
By Corollary \ref{X is non-essential}, it follows that $n$ is even.
\end{proof}

\begin{lemma}
\label{verhnij cykl}
Let $\alpha\in\BCH(X^2)$ be a cycle containing
$
\beta=h^{\jw_{q-1}}\times l_{\jw_q-1}
$
for some $q\in S=\{1,2,\dots,\Hight\}$
(this $\beta$ is the top of the $q$-th left shell triangle).
Then
$\alpha$ also contains the transposition of $\beta$.
\end{lemma}

\begin{proof}
Replacing $F$ by the field $F_{q-1}$ of the generic splitting
tower of $F$, $X$ by $X_{q-1}$, and $\alpha$ by $\pr^2(\alpha_{F_{q-1}})$,
we come to the situation where $q=1$.

Then we replace $\alpha$ by its homogeneous component containing
$\beta$ and apply to it Lemma \ref{dyrki v verhnem cykle} (with $k=\iw_1$).
Let us assume that the transposition of $\beta$ is
not contained in $\alpha$.

By Lemma \ref{dyrki v verhnem cykle}
$\alpha$ does not contain any of the essential basis cycles having $h^i$ with
$0<i<\iw_1$ as a factor;
therefore
the number of the essential basis elements contained in $\alpha$ and
the number of the essential basis elements contained in
$\pr^2(\alpha_{F(X)})\in\BCH(X_1^2)$
differ by $1$.
In particular, these two numbers
have
different parity.
However, the number of the essential basis elements contained in $\alpha$ is even
by Lemma \ref{chetnost' chisla tochek}.
By the same lemma, the number of the essential basis elements contained in
$\pr^2(\alpha_{F(X)})$ is even too.
\end{proof}

\begin{dfn}
\label{def-primordial}
A minimal cycle $\alpha\in\BCH^{\leq D}(X^2)$ is called {\em
primordial}, if it is not a positive order derivative of another
rational cycle.
\end{dfn}

\begin{lemma}
\label{symmetry}
Let $\alpha\in\BCH(X^2)$ be a minimal cycle.
Assume that for some $q\in S$, the cycle $\alpha$ contains
$h^{\jw_{q-1}}\times l_{\jw_q-1}$.
Then $\alpha$ is symmetric and primordial.
\end{lemma}

\begin{proof}
The cycle
$
\alpha\cap t(\alpha)
$
(where $t(\alpha)$ is the transposition of $\alpha$;
intersection of cycles is defined in Lemma \ref{intersection})
is symmetric, rational (Lemma \ref{intersection}),
contained in $\alpha$, and, by Lemma \ref{verhnij cykl},
still contains $h^{\jw_{q-1}}\times
l_{\jw_q-1}$
(in particular, $\alpha\cap t(\alpha)\ne0$).
It coincides with $\alpha$ by the minimality of $\alpha$.

It is easy to ``see'' that $\alpha$ is primordial looking at the
picture of Remark \ref{big piramid} (because $\alpha$ contains the top point of
some shell triangle).
Nevertheless, let us do the proof by formulae.
If there exists a rational cycle $\beta\ne\alpha$ such that $\alpha$ is a
derivative of $\beta$, then there exists a rational cycle $\beta'$
such that $\alpha$ is an order one derivative of $\beta'$, that
is, $\alpha=\beta'\cdot(h^0\times h^1)$ or $\alpha=\beta'\cdot(h^1\times
h^0)$.
In the first case $\beta'$ should contain the basis cycle
$h^{\jw_{q-1}}\times l_{\jw_q}$,
while in the second case $\beta'$ contains
$h^{\jw_{q-1}-1}\times l_{\jw_q-1}$.
However, these both cases are not possible by Lemma \ref{dyrki v verhnem cykle}
(take $k=\iw_q+1$ with $i=0$ for the first case and $i=\iw_{q-1}-1$ for the second
case).
\end{proof}

It is easy to see that a cycle $\alpha$ with the property of Lemma
\ref{symmetry} exists at least for $q=1$:

\begin{lemma}
\label{vneshnij verhnij cykl}
There exists a cycle in $\BCH_{D+\iw_1-1}(X^2)$ containing $h^0\times
l_{\iw_1-1}$.
\end{lemma}

\begin{proof}
Take a preimage of $l_{\iw_1-1}\in\BCH(X_{F(X)})$
under the surjection
$\BCH(X^2)\onto\BCH(X_{F(X)})$ given by the pull-back with respect
to the morphism $X_{F(X)}\to X^2$ produced by the generic point of
the first factor of $X^2$.
\end{proof}

The following lemma is proved already in \cite{i_1} (under the name of ``Vishik's
principle''), but only for odd-dimensional quadrics and by a
different as here method.

\begin{lemma}
\label{pairs}
For any cycle $\rho\in\BCH_\D(X^2)$, any $q\in S$, and any
$i\in[1,\;\iw_q]$,
the element
$
h^{\jw_{q-1}+i-1}\times
l_{\jw_{q-1}+i-1}
$
is contained in $\rho$ if and only if
the element
$
l_{\jw_q-i}\times
h_{\jw_q-i}
$
is contained in $\rho$.
\end{lemma}

\begin{proof}
Clearly, it is enough to prove Lemma \ref{pairs} for $q=1$ only.
Let $\rho$ be a counter-example to Lemma \ref{pairs}.
Replacing $\rho$ by its transposition, if necessary, we can come to the
situation where
$\rho\ni l_{\iw_1-i}\times h^{\iw_1-i}$ and
$\rho\not\ni h^{i-1}\times l_{i-1}$ for some $i\in[1,\;\iw_1]$.

Let $\alpha$ be the cycle of Lemma \ref{vneshnij verhnij cykl}.
The composition
$
\beta=\big(t(\alpha)\cdot(h^{i-1}\times h^0)\big)\compose\rho
$
(where $t(\alpha)$ is the transposition of $\alpha$)
is a rational homogeneous cycle
containing $l_{\iw_1-i}\times h^0$ and
not containing $h^{i-1}\times l_{\iw_1-1}$.
Therefore
$
\gamma=t\big(\beta\cdot(h^{\iw_1-i}\times h^0)\big)
$
is a rational
cycle
containing $h^0\times l_0$ and
not containing
$l_{\iw_1-1}\times h^{\iw_1-1}$.
It follows that the composition
$\alpha\compose\gamma$ is a rational cycle containing
$h^0\times l_{\iw_1-1}$ and not containing $l_{\iw_1-1}\times h^0$,
a contradiction with Lemma \ref{verhnij cykl}.
\end{proof}

To announce the result which follow,
we prefer to use the language of picture rather then the language
of formulae:

\begin{cor}
\label{more pairs}
The diagram of an arbitrary $\alpha\in\BCH^{\leq D}(X^2)$ has the
following property:
for any $q\in S$ and any integers $i\geq1$ and $k\geq0$,
the $i$-th point of the $k$-th row of the $q$-th {\em left} shell
triangle is marked if and only if the $i$-th point of the $k$-th row of the $q$-th
{\em right} shell triangle is marked
(see Definition \ref{def-shell triangles} for the agreement on counting the rows
and the points of the shell triangles).
\end{cor}

\begin{proof}
The case of $k=0$ is treated in Lemma \ref{pairs}
(while Lemma \ref{verhnij cykl} treats the case of ``maximal''
$k$).
The case of an arbitrary $k$ is reduced to the case of $k=0$ by
taking a $k$-th order derivative of $\alpha$.
\end{proof}

\begin{rem}
\label{half-diagram}
By Corollary \ref{more pairs},
it follows that the diagram of a cycle in $\BCH^{\leq D}(X^2)$
is determined by, say, the left half of itself.
\end{rem}

\begin{example}
\label{new proof of i1}
As an application of the results on $X^2$ obtained by now (first
of all, of Corollary \ref{more pairs}), we give a short (simpler
as the original) proof of the main result of \cite{i_1}, which can
be stated as follows:
if $\phi$ is an anisotropic quadratic form and $2^r$ is the biggest
power of $2$ dividing the difference $\dim(\phi)-\iw_1(\phi)$,
then $\iw_1(\phi)\leq2^r$.
For the proof, assume that $\iw_1=\iw_1(\phi)>2^r$ and consider
the Steenrod operation $S^{2^r}(\alpha)$ of a homogeneous cycle
$\alpha\in\BCH^{\leq D}(X^2)$ containing $h^0\times l_{\iw_1-1}$
(for the existence of $\alpha$ see Lemma \ref{vneshnij verhnij
cykl};
note that $S^{2^r}(\alpha)$ is still inside of $\BCH^{\leq D}(X^2)$ just because
of the inequality $\iw_1>2^r$).
Since
$$
S^{2^r}(h^0\times l_{\iw_1-1})=h^0\times S^{2^r}(l_{\iw_1-1})=
\binom{\dim(\phi)-\iw_1}{2^r}\cdot(h^0\times l_{\iw_1-1-2^r})
$$
and the binomial coefficient is odd, we get that
$S^{2^r}(\alpha)\ni h^0\times l_{\iw_1-1-2^r}$.
On the other hand,
$\alpha\not\ni l_{\iw_1-1+i}\times h^i$ for any $i\in[1,\iw_1-1]$
by Lemma \ref{dyrki v verhnem cykle};
consequently, $S^{2^r}(\alpha)\not\ni l_{\iw_1-1-2^r+i}\times h^i$
for these $i$;
in particular, this is so for $i=\iw_1-1$.
Now, applying Corollary \ref{more pairs}, we get that
$S^{2^r}(\alpha)\not\ni h^0\times l_{\iw_1-1-2^r}$, a
contradiction.
\end{example}

The following lemma generalizes Lemma \ref{vneshnij verhnij cykl}:

\begin{lemma}
\label{suschestvovanie verhnego cykla}
Let $q\in S$.
Assume that the group $\BCH_D(X^2)$ contains a cycle $\gamma$ such
that
\begin{enumerate}
\item
$\gamma$ does not contain any $h^i\times l_i$ with $i<
\jw_{q-1}$;
\item
$\gamma$ contains $h^i\times l_i$ for some
integer $i\in[\jw_{q-1},\;\jw_q)$
(note that the interval is semi-open).
\end{enumerate}
Then the group $\BCH_{D+\iw_q-1}(X^2)$ contains a cycle $\alpha$
such that
$\alpha\ni h^{\jw_{q-1}}\times
l_{\jw_q-1}$ and
$\alpha\not\ni h^i\times l_{i+\iw_q-1}$ for any $i<\jw_{q-1}$.
\end{lemma}

\begin{proof}
We use an induction on $q$.
In the case of $q=1$, the assumption of Lemma \ref{suschestvovanie verhnego cykla}
is always true (think of $\gamma=\Delta$);
the cycle $\alpha$ is constructed in Lemma \ref{vneshnij verhnij
cykl}.
In the remaining part of the proof we assume that $q>1$.

Let $i$ be the smallest integer such that
$\gamma\ni h^{\jw_{q-1}+i}\times
l_{\jw_{q-1}+i}$.
As a first step, we proof that the group $\BCH^{\leq D}(X^2)$ contains
a cycle $\alpha'$ containing
$h^{\jw_{q-1}+i}\times l_{\jw_q-1}$
and none of $h^j\times l_?$ with $j<\jw_{q-1}+i$
(if $i=0$ then we can take $\alpha=\alpha'$ and finish the proof).

Applying the induction hypothesis to the quadric $X_1$
with the cycle $\pr^2(\gamma_{F(X)})\in\BCH(X_1^2)$
(and using the inclusion
homomorphism of Corollary \ref{old isotropic}), we get a
cycle in $\BCH_{D+\iw_q-1}(X_{F(X)}^2)$ containing
$
h^{\jw_{q-1}}\times l_{\jw_q-1}\;.
$
One of its derivatives is a homogeneous cycle in $\BCH(X^2_{F(X)})$
containing
$
h^{\jw_{q-1}+i}\times l_{\jw_q-1}\;.
$
Therefore the group $\BCH(X^3)$ contains a homogeneous cycle
containing
$
h^0\times h^{\jw_{q-1}+i}\times l_{\jw_q-1}\;.
$
Considering it as a correspondence of the middle factor of
$X^3$ into the product of two outer factors, composing it with
$\gamma$, and taking the pull-back with respect to the first diagonal
$X^2\to X^3$, we get the required cycle $\alpha'$.

The highest order derivative $\alpha'\cdot(h^{\iw_q-1-i}\times h^0)$ of $\alpha'$
contains
$
h^{\jw_q-1}\times l_{\jw_q-1}\;.
$
By Lemma \ref{pairs}, it also contains
$
l_{\jw_{q-1}}\times h^{\jw_{q-1}}\;.
$
Therefore its transposition contains
$
h^{\jw_{q-1}}\times l_{\jw_{q-1}}\;.
$
Replacing $\gamma$ by the constructed rational cycle, we come to
the situation with $i=0$
(see the second paragraph of the proof),
finishing the proof.
\end{proof}

We come to the main result of on the structure of $\BCH^{\leq
D}(X^2)$ for an arbitrary anisotropic projective quadric $X$:

\begin{thm}
\label{kruto!}
The set of the primordial (see Definition \ref{def-primordial}) cycles
$\Pi\subset\BCHE^{\leq D}(X^2)$ has the following properties.
\begin{enumerate}
\item
All derivatives of all cycles of $\Pi$ are minimal and pairwise
different;
they form a basis of $\BCHE^{\leq D}(X^2)$.
\item
Every cycle in $\Pi$ is symmetric.
\item
For every $\pi\in\Pi$, there exists one and only one
$q=f(\pi)\in S=\{1,2,\dots,\Hight\}$ such
that
\begin{enumerate}
\item
$\dim(\pi)=D+\iw_q-1$;
\item
$\pi\not\ni h^i\times l_{i+\iw_q-1}$ for any
$i<\jw_{q-1}$;
\item
$\pi\ni h^{\jw_{q-1}}\times
l_{\jw_q-1}$.
\end{enumerate}
\item
The map $f\!:\Pi\to S$ thus obtained is injective, its
image consists of $q\in S$ such that there exists a cycle
$\alpha\in\BCH^{\leq D}(X^2)$ satisfying
$\alpha\ni h^i\times l_?$ for some $i\in[\jw_{q-1},\;\jw_q)$ and
$\alpha\not\ni h^i\times l_?$ for any $i\in[0,\;\jw_{q-1})$
(in particular, $f(\Pi)\ni1$).
\end{enumerate}
\end{thm}

%

\begin{proof}
We construct
a chain of subsets
$$
\emptyset=\Pi_0\subset \Pi_1\subset\dots\subset \Pi_\Hight
$$
of the set $\Pi$ such that
for every $q\in S$,
all highest derivatives of all cycles of $\Pi_q$ are minimal and pairwise
different, and their sum contains $h^i\times l_i$ for all $i<\jw_q$.
The procedure looks as follows.
If for some $q\in S$ the set $\Pi_{q-1}$ is
already constructed, we decide whether we
set $\Pi_q=\Pi_{q-1}\cup\{\pi\}$ with certain cycle $\pi$
or we set $\Pi_q=\Pi_{q-1}$.
To make this decision, we consider the sum $\alpha$ of all highest derivatives of all
cycles of $\Pi_{q-1}$.
We know that $\alpha$ contains $h^i\times l_i$ for all
$i\in[0,\;\jw_{q-1})$.
If $\alpha$ also contains
$h^i\times l_i$ for all
$i\in[\jw_{q-1},\;\jw_q)$, then we set $\Pi_q=\Pi_{q-1}$;
otherwise the cycle $\gamma=\alpha+\Delta$
satisfies the hypothesis of Lemma \ref{suschestvovanie verhnego
cykla}, and we set $\Pi_q=\Pi_{q-1}\cup\{\pi\}$ with $\pi$ being
the minimal cycle containing
$
h^{\jw_{q-1}}\times l_{\jw_q-1}
$
($\pi$ exists and has Property (3b) by Lemma \ref {suschestvovanie verhnego cykla};
$\pi$ is primordial by Lemma \ref{symmetry}).

The set $\Pi_\Hight$ thus constructed has all properties claimed for
$\Pi$ in Theorem \ref{kruto!}.
Indeed the elements of $\Pi_\Hight$ are symmetric by
Lemma \ref{symmetry}.
%
The sum of all highest derivatives of all elements of $\Pi_\Hight$ contains
$h^i\times l_i$ for all $i$; therefore this sum also contains
the remaining basis elements
$l_i\times
h^i$ for all $i$ (see Lemma \ref{pairs}).
It follows that every $D$-dimensional minimal cycle is a
derivative of an element of $\Pi_\Hight$.
Consequently, every minimal cycle in $\BCH^{\leq D}(X^2)$
is a derivative of of a cycle of $\Pi_\Hight$.
It follows that $\Pi_\Hight=\Pi$.
All minimal cycles form a basis according to
Proposition \ref{minimal form basis}.
\end{proof}

As easy as important information on relations between the
primordial cycles on $X^2$ and on $X_1^2$ is as follows:

\begin{prop}
\label{neravenstva}
Let $\Pi$ be the set of all primordial cycles for $X$;
let $\Pi_1$ be the set of all primordial cycles for $X_1$.
As usual we set $\iw_1=\iw_1(X)$.
One has:
\begin{enumerate}
\item
$\#\Pi-1\leq\#\Pi_1$;
\item
if $\Pi\not\ni h^0\times l_{\iw_1-1}+ l_{\iw_1-1}\times h^0$, then
$\#\Pi\leq\#\Pi_1$.
\end{enumerate}
\end{prop}

\begin{proof}
Let us extend the function $f\!:\Pi\to S$ on
the set of all non-zero cycles in $\BCHE^{\leq D}(X^2)$,
defining $f(\alpha)$
as the minimal $q\in S$ such that $\alpha\ni h^i\times l_?$ for
some $i\in[\jw_{q-1},\;\jw_q)$ and $\alpha\not\ni h^i\times l_?$
for any $i\in[0,\;\jw_{q-1})$.
By Item 4 of Theorem \ref{kruto!} (which is a consequence of Lemma
\ref{suschestvovanie verhnego cykla}), the image of the extended
$f$ coincides with $f(\Pi)$.
Let $f_1\!:\BCHE^{\leq D}(X_1^2)\to S_1$ be the same map for the
quadric $X_1$.
We denote as $\Pi'$ the set $\Pi$ without the primordial cycle
containing $h^0\times l_{\iw_1-1}$
(this is the primordial cycle whose image under $f$ is $1$).
For any $\pi\in\Pi'$
the cycle
$\pr^2(\pi)\in\BCHE(X_1^2)$ is non-zero and
$f_1(\pr^2(\pi))=f(\pi)-1$.
It follows that $\#\Pi_1=\# f_1(\Pi_1)=\# \Im(f_1)\geq \#
f(\Pi')=\#\Pi'=\#\Pi-1$,
and the first statement of Proposition \ref{neravenstva} is
proved.


If now $\Pi\not\ni h^0\times l_{\iw_1-1}+l_{\iw_1-1}\times h^0$, then
$\pr^2(\pi_{F(X)})$ is non-zero for {\em every} $\pi\in\Pi$.
Note that for the cycle $\pi\in\Pi$ containing $h^0\times
l_{\iw_1-1}$, one has $f_1(\pr^2(\pi))\not\in f_1(\pr^2(\Pi'))$
(because $\pi$ is disjoint with all derivatives of the cycles of $\Pi'$ and,
consequently, $\pr^2(\pi)$ is disjoint with all derivatives of the cycles of
$\pr^2(\Pi')$).
Therefore $\#\Pi\leq\#\Pi_1$,
and the second statement of Proposition \ref{neravenstva} is
proved as well.
\end{proof}

We need some more notation.

\begin{dfn}
For any $r\geq1$,
the symmetric group $S_r$ acts on the group $\CH(\bar{X}^r)$
by permutations of factors of $\bar{X}^r$.
If $\alpha\in\CH(\bar{X}^r)$, we write $\Sym(\alpha)$ for the ``symmetrization''
of $\alpha$, that is,
$$
\Sym(\alpha)=\Sum_{s\in S_r} s(\alpha)\;.
$$
\end{dfn}

\begin{dfn}
\label{def-small}
A non-zero anisotropic quadratic form $\phi$ over $F$ is said to be {\em small} if for some
positive integer $n$ (which is uniquely determined by $\dim(\phi)$
by the Arason-Pfister theorem)
one has $\phi\in I^n$ while
$\dim\phi<2^{n+1}$.
A projective quadric is {\em small} if so is the
corresponding quadratic form.
\end{dfn}

The following result is an extended version of \cite[thm. 4.1]{Vishik-In}.

\begin{prop}
\label{known}
Let $X$ be a small $2d$-dimensional quadric
of the first Witt index $a=\iw_1(X)$.
Then
\begin{enumerate}
\item
the integer $a$ divides all the higher Witt indices $\iw_1,\dots,\iw_\Hight$ of $X$;
in particular, it divides $d+1=\iw_1+\dots+\iw_\Hight$;
\item
the cycle
$$
\pi=
\Sym\left(\Sum_{i=1}^{(d+1)/a}h^{(i-1)a}\times l_{ia-1}\right)
\in\CH_{2d+a-1}(\bar{X}^2)
$$
is rational;
\item
\label{item 3}
moreover,
for every
$k\geq0$
the Chow group
$\BCH_{2d+k}(X^2)$ is generated by the non-essential basis elements and
the cycles
$$
\pi\cdot(h^{j-1}\times h^{a-k-j})\;,\;\;j=1,2,\dots,a-k
$$
(in particular, for $k\geq a$, this Chow group consists of the non-essential
elements only).
\end{enumerate}
\end{prop}

\begin{proof}
Let $\Pi$ be the set of primordial cycles.
It suffices to show that $\#\Pi=1$ (then the unique element
$\pi\in \Pi$ automatically has the form and the property
required).
We prove it using an induction on $\Hight=\Hight(X)$.
If  $\Hight=1$, then $\#\Pi=1$,
since generally $1\leq \#\Pi\leq \Hight$.

Now we assume that $\Hight\geq2$.
Let us consider the quadric $X_1$ (over the field $F(X)$)
and let $\Pi_1$ be the set of primordial cycles for $X_1$.
Then
$\#\Pi_1=1$ by the induction hypothesis, and we get what we need
by Item 2 of Proposition \ref{neravenstva},
if we check that
the cycle $\Sym(h^0\times l_{a-1})$ is not rational.
By Item \ref{binary size} of Proposition \ref{list},
this cycle can be rational only if
the integer $2d-(a-1)+1=2d-a+2$ is a power of $2$.
Since however
$$
2^n\leq
\dim(\phi_1)=2d+2-2a<\underline{\underline{2d-a+2}}<2d+2=\dim(\phi)\leq2^{n+1}\;,
$$
the integer $2d-a+2$ is {\em not} a power of $2$.
\end{proof}

\begin{rem}
Proposition \ref{known} holds also for anisotropic $\phi$ with $[\phi]\in
I^n$ and $\dim(\phi)=2^{n+1}$ (the same proof is valid for such $\phi$ as well).
\end{rem}

\section
{Cycles on $X^3$}

Let $\phi$ be a small quadratic form and let $n$ be the positive integer such that
$[\phi]\in I^n$ while $\dim(\phi)<2^{n+1}$.
We recall that $X$ stays for the projective quadric given by $\phi$.
Let us write down the dimension of $\phi$ as a sum of powers of $2$:
$$
\dim(\phi)=2^n+2^{n_1}+\dots+2^{n_m}\;,\;\;
n>n_1>\dots>n_m\geq1\;.
$$
In this section we {\em assume} that $m\geq2$, that the height of $\phi$ is at least $3$,
and that
the first two higher Witt indices of $\phi$ are
as follows:
$\iw_1(\phi)=2^{n_m-1}$ and
$\iw_2(\phi)=2^{n_{m-1}-1}$.
To simplify the formulae which follow, we introduce the notation
$$
a=\iw_1(\phi);\;\;
b=\iw_2(\phi);\;\;
c=\iw_3(\phi)
$$
and
$$
d=\dim(X)/2=2^{n-1}+2^{n_1-1}+\dots+2^{n_m-1}-1\;.
$$

Here is our main construction:

\begin{prop}[{cf. \cite[prop. 2.7]{third proof}}]
\label{mu}
The group $\BCH(X^3)$ contains a homogeneous cycle
$$
\mu=
\Sym\left(\Sum_{i=1}^{(d-a+1)/b}h^0\times h^{(i-1)b+a}\times l_{ib+a-1}\right)
+\mu'
\;,
$$
where $\mu'$ is a sum of only those essential
basis elements which
contain neither $h^0$ nor $h^i$ with $a\not|$ $i$.
\end{prop}

\begin{proof}
Let $X_1$ be the projective quadric (over the field $F(X)$) given by the anisotropic part
of the form $\phi_{F(X)}$.
Applying Item 2 of Proposition \ref{known} to the quadric $X_1$ (taking in account that
$\iw_1(X_1)=\iw_2(X)=b$ and $\dim(X_1)=2(d-a)$), we , in particular, get that the group
$\BCH(X_1^2)$ contains the cycle
$$
\beta'=
\Sym\left(\Sum_{i=1}^{(d-a+1)/b}h^{(i-1)b}\times l_{ib-1}\right)
\;.
$$
Therefore (see Corollary \ref{old isotropic}),
the group $\BCH(X_{F(X)}^2)$ contains the cycle
$$
\beta=\inc^2(\beta')=
\Sym\left(\Sum_{i=1}^{(d-a+1)/b}h^{(i-1)b+a}\times l_{ib+a-1}\right)
\;.
$$
The pull-back homomorphism $g_1^*\!:\BCH(X^3)\to\BCH(X_{F(X)}^2)$
with respect to the morphism
$g_1\!:X^2_{F(X)}\to X^3$,
given by the generic point of
the first factor of $X^3$, is surjective.
Therefore, there exists a homogeneous cycle $\mu\in\BCH(X^3)$ such that
$g_1^*(\mu)=\beta$.
Note that $g_1^*$ sends every basis cycle of the type
$h^0\times\zeta\times\xi$ to $\zeta\times\xi$ while killing the other basis
elements.
Consequently we have
$$
\mu=h^0\times
\Sym\left(\Sum_{i=1}^{(d-a+1)/b}h^{(i-1)b+a}\times l_{ib+a-1}\right)
+ \epsilon\;,
$$
where $\epsilon$
is a sum of some basis cycles which do not have $h^0$ on the first factor place.

We now proceed by transforming the cycle $\mu$ in such
a way that $\mu$ does not leave the group
$\BCH(X^3)$ and $g_1^*(\mu)$ remains the same.

By Proposition \ref{known} (now applied to $X$ itself), the cycle
\begin{multline*}
\gamma=
\Sym\left(\Sum_{i=1}^{(d+1)/a}h^{(i-1)a}\times l_{ia-1}\right)
\cdot(h^0\times h^{a-1})=\\
\Sum_{i=1}^{(d+1)/a}(h^{(i-1)a}\times l_{(i-1)a}+l_{ia-1}\times h^{ia-1})
\end{multline*}
is in $\BCH(X^2)$.
Considering it as a correspondence, we replace $\mu$ by the composition
$\mu\compose\gamma$,
where $\mu\in\BCH(X_1\times X_2\times X_3)$
is considered as a correspondence from $X_1$ to
$X_2\times X_3$ (all $X_i$ are copies of $X$).
Now a basis element  $h^i\times?\times?$ occurs in the decomposition of $\mu$ only if
$i$ is divisible by $a$ (see Lemma \ref{comp-corr}).

Considering $\mu$ as a correspondence from $X_2$ to $X_1\times X_3$ and replacing it by
the composition $\mu\compose\gamma$,~\footnote{Strictly speaking, this is
$t_{12}\big(t_{12}(\mu)\compose\gamma\big)$,
where $t_{12}$ is the automorphism of $\BCH(X^3)$ induced by the transposition of
the first two factors of $X^3$.}
we come to the situation where
a basis element  $?\times h^i\times?$ occurs in the decomposition of $\mu$ only if
$i$ is divisible by $a$ (while all previously established properties of $\mu$ still hold).

Finally, considering $\mu$ as a correspondence from $X_3$ to $X_1\times X_2$ and replacing it by
the composition $\mu\compose\gamma$, we come to the situation where
a basis element  $?\times ?\times h^i$ occurs in the decomposition of $\mu$ only if
$i$ is divisible by $a$ (while all the previously established properties of $\mu$ still hold).

The last change we apply to $\mu$ is as follows:
we remove all non-essential basis cycles
in the decomposition of $\mu$.
We claim that now our cycle $\mu$ has the required shape.

Let us write $\mu_0$ for the sum of those summands in the decomposition of $\mu$ which
have $h^0$ as at least one factor.
To finish the proof of the proposition, it suffices to check that
$$
\mu_0=
\Sym\left(\Sum_{i=1}^{(d-a+1)/b}h^0\times h^{(i-1)b+a}\times l_{ib+a-1}\right)
\;.
$$
First of all let us check that none of the 3 basis cycles
obtained from $h^0\times h^0\times l_{b-1}$ by a permutation of factors appears
in the decomposition of $\mu_0$.
We assume that the cycle $h^0\times h^0\times l_{b-1}$ does appear and
we pull-back $\mu$ with respect to the morphism
$g_{12}\!:X_{F(X\times X)}\to X^3$ given by the generic point of the product of the
first two factors of $X^3$.
We get
$$
\BCH(X_{F(X\times X)})\ni g_{12}^*(\mu)= g_{12}^*(h^0\times
h^0\times l_{b-1})= l_{b-1}
$$
showing that the Witt index of the quadric
$X_{F(X\times X)}$ is at least $b$
(see Corollary \ref{determining witt index}).
However this Witt index coincides with
$\iw_1(X)=a$ and $a$ is smaller than $b$ (actually $a\leq b/2$).
The contradiction obtained shows that the cycle $h^0\times h^0\times l_{b-1}$ does
not appear
in the decomposition of $\mu_0$.
For every permutation of $h^0\times h^0\times l_{b-1}$ we get the same result simply by
changing in the appropriate way the numeration of the factors.

It follows that $\mu_0=\mu_1+\mu_2+\mu_3$, where $\mu_i$ is the sum of summands in the
decomposition of $\mu$ such that $h^0$ is their $i$-th factor.
By the construction of $\mu$ we known that
$$
\mu_1=
h^0\times
\Sym\left(\Sum_{i=1}^{(d-a+1)/b}h^{(i-1)b+a}\times l_{ib+a-1}\right)\;,
$$
and it suffices to check that $\mu_2=t_{12}(\mu_1)$ and $\mu_3=t_{13}(\mu_1)$
with $t_{1i}$ staying for the
automorphism of the Chow group $\CH(\bar{X}^3)$ given by the
transposition of the first and $i$-th factor of $\bar{X}^3$.

In order to see that $\mu_2=t_{12}(\mu_1)$, we pull-back the cycle $\mu$ to
$X^2$ with respect to the morphism
$$
\delta_1\!:X^2\to X^3\;,\;\;x_1\times x_2\mapsto x_1\times
x_1\times x_2
$$
given by the diagonal map
of the first factor of $X^2$ into the product of the first two factors of $X^3$.
The decomposition of the homogeneous cycle
$\delta_1^*(\mu)\in\BCH_{2d+b-1}(X^2)$ does not contain
any non-essential cycle.
Therefore, since $b>a$, $\delta_1^*(\mu)=0$ by Proposition \ref{known}.
On the other hand, $\delta_1^*(\mu_1)$ contains $h^a\times l_{b+a-1}$ while
neither $\delta_1^*(\mu_3)$ nor $\delta_1^*(\mu-\mu_0)$ do.
It follows that $\delta_1^*(\mu_2)$ contains $h^a\times l_{b+a-1}$ as well
and consequently $\mu_2$ contains the basis cycle
$h^a\times h^0\times l_{b+a-1}$.
Now we use the pull-back with respect to the morphism
$g_2\!:X^2_{F(X)}\to X^3$ given by the generic point of the second factor of $X^3$.
The homogeneous cycle
$g_2^*(\mu)=g_2^*(\mu_2)$ lies in $\BCH(X^2_{F(X)})$, contains
the basis cycle $h^a\times l_{b+a-1}$, and does not contain
any non-essential basis element.
Passing to the anisotropic part $X_1$ of $X_{F(X)}$
and using Corollary \ref{old isotropic}, we get a
homogeneous cycle $\eta$
in $\BCH(X_1^2)$, namely $\eta=\pr^2(g_2^*(\mu))$,
which contains $h^0\times l_{b-1}$ and does not contain any
non-essential cycle.
Note that $g_2^*(\mu_2)$ is in the image of
$\inc^2\!:\CH(\bar{X}_1^2)\to\CH(\bar{X}_{F(X)}^2)$,
so that $\mu_2$ can be reconstructed from $\eta$.

By Proposition \ref{known} it follows that
$$
\eta=
\Sym\left(\Sum_{i=1}^{(d-a+1)/b}h^{(i-1)b}\times l_{ib-1}\right)
\;.
$$
Consequently
$$
g_2^*(\mu_2)=\inc^2(\eta)=
\Sym\left(\Sum_{i=1}^{(d-a+1)/b}h^{(i-1)b+a}\times l_{ib+a-1}\right)
$$
and $\mu_2=t_{12}(\mu_1)$.

The equality $\mu_3=t_{13}(\mu_1)$ is checked similarly.
\end{proof}

We remark that the ``defect part'' $\mu'$ of the cycle $\mu$ does not appear
in \cite[prop. 2.7]{third proof} when working with a small quadric of height $2$.
In our case here, the height of $X$ is at least $3$, $\mu'$ does really exist and
represents an additional difficulty.
The main observation which is crucial to overcome this difficulty is as follows:

\begin{lemma}
\label{mu'}
Let $\mu'$ be
as in Proposition \ref{mu}.
In the decomposition of $\mu'$ we consider the basis elements with $h^a$
on the $i$-th factor place and write $\mu'_i$ for their sum.
Then each of the cycles
$\mu'_1$, $t_{12}(\mu'_2)$, and $t_{13}(\mu'_3)$ is the sum of
\underline{some} of the following $(c-b)/a$ elements
$$
\chi_j=
h^a\times\left(
\Sym\left(
\Sum_{i=1}^{(d-b-a+1)/c}h^{(i-1)c+b+a}\times l_{ic+b+a-1}
\right)\cdot\Big(h^{(j-1)a}\times h^{c-b-ja}\Big)
\right)
$$
where $j\in\{1,2,\dots,(c-b)/a\}$.
\end{lemma}

\begin{proof}
Clearly, it suffices to prove the statement on $\mu'_1$
(the statements on $\mu'_2$ and $\mu'_3$ are proved in the same way
interchanging the roles of the three factors of $X^3$).

Let us go over the function field $F(X)$.
We still have $\mu\in\BCH(X^3_{F(X)})$.
Therefore
$\pr^3_X(\mu)\in\BCH(X_1^3)$, where $X_1$ is the anisotropic part of $X_{F(X)}$
and $\pr^3_X\!:\CH(\bar{X}_{F(X)}^3)\to\CH(\bar{X}_1^3)$ is the projection
of Corollary \ref{old isotropic}.
We note that $\pr^3_X(\mu)=\pr^3_X(\mu')$.
Moreover, $\mu'$ is in the image of the inclusion
$\inc^3_X\!:\CH(\bar{X}_1^3)\to\CH(\bar{X}_{F(X)}^3)$
because every $h^i$ which is a factor of a basis element in the
decomposition of $\mu'$ has $i\geq a$ and every $l_i$
which is a factor of a basis element in the
decomposition of $\mu'$ has $i\geq a$ as well (just look at the dimension
of $\mu'$).
Therefore $\mu'$ can be reconstructed from its image under $\pr^3_X$,
namely, $\mu'=\inc^3_X(\pr^3_X(\mu'))$.

Now we move from $\BCH(X_1^3)$ to $\BCH\big((X_1)_{F(X)(X_1)}^2\big)$ using
$g_1^*$ (the pull-back with respect to the morphism
given by the generic point of the first factor of $X_1^3$).
Note that $g_1^*(\pr^3_X(\mu'))=g_1^*(\pr^3_X(\mu'_1))$
and the cycle $\pr^3_X(\mu'_1)$ can be reconstructed from its image under
$g_1^*$.
Moreover the cycle $g_1^*(\pr^3_X(\mu'_1))$ is in the image of the inclusion
$\inc^2_{X_1}\!:\CH(\bar{X}_2^2)\to\CH\big((\bar{X}_1)_{F(X_1)}^2\big)$,
where $X_2$ is the anisotropic part of $(X_1)_{F(X)(X_1)}$.
In order to see it, we note that every basis cycle in the decomposition of
$g_1^*(\pr^3_X(\mu'_1))$ is equal to $h^{(i-1)a}\times l_{b+ia-1}$
for some $i\geq1$.
Clearly, such a basis cycle is
in the image of $\inc^2_{X_1}$ if and only if $(i-1)a\geq b$.
So, if the cycle $g_1^*(\pr^3_X(\mu'_1))$ is not in the image of $\inc^2_{X_1}$,
then its decomposition contains the basis element
$h^{(i-1)a}\times l_{b+ia-1}$ with some $i$ such that $(i-1)a<b$.
It follows that the decomposition of the rational cycle
$$
\pr_{2d-2a-(i-1)a,a}(g_1^*(\pr^3_X(\mu'_1)))\in\BCH(X_2)
$$
contains
$l_{ia-1}$.
This is a contradiction because $X_2$ is anisotropic and therefore
the group $\BCH(X_2)$
does not contains essential elements (Corollary \ref{X is
non-essential}).

So, one can reconstruct the cycle $g_1^*(\pr^3_X(\mu'_1))$ from
its image under the projection
$\pr^2_{X_1}\!:\BCH\big((X_1)_{F(X)(X_1)}^2\big)\to\BCH(X_2^2)$
and for our purposes
it is sufficient to determine this
image.
To do so, we apply Proposition \ref{known}
to the quadric $X_2$ getting that the cycle
$$
(\pr^2_{X_1}\compose g_1^*\compose \pr^3_X)(\mu'_1)\in\BCH_{2d-b-a-1}(X_2^2)
$$
is the sum of some essential
generators of the group $\BCH_{2d-b-a-1}(X_2^2)$ indicated in Item 3 of
Proposition \ref{known}
(we note that $\dim(X_2)=2(d-b-a)$ so that $2d-b-a-1=\dim(X_2) +b+a-1\leq\dim(X_2)+c-1$).
Finally, taking in account that $h^i$ for a given $i$ can be a factor of a basis
element appearing in the decomposition of
$(\pr^2_{X_1}\compose g_1^*\compose \pr^3_X)(\mu'_1)$
only if $i$ is divisible by $a$,
we get the desired description of the cycle $\mu'_1$.
\end{proof}

\section
{Proof of Conjecture \ref{conj}}
\label{Proof of Conjecture}

In this {\S} we prove Conjecture \ref{conj}.
Suppose that this conjecture is not true, that is, over some field $F$
and for some positive integer $n$, there
exists a quadratic form $\phi$ over $F$ with
$[\phi]\in I^n$  and with $\dim(\phi)$ prohibited
by Conjecture \ref{conj}.
Note that $n$ is at least $4$ (see \S\ref{Introduction}).
In the splitting pattern of the form $\phi$,
let us choose the smallest number $\dim(\phi_E)_{\an}$
prohibited by Conjecture \ref{conj}.
Let us replace the form $\phi$ by this $(\phi_E)_{\an}$ (and $F$ by $E$)
and write $X$ for the projective quadric given by the new $\phi$.
Note that
$\dim(\phi)>2^n+2^{n-1}$
(\cite{Vishik-Lens} (the original proof), or \cite[thm. 4.4]{KM},
or \cite{third proof}).
Moreover,
$$
\dim(\phi_{F(X)})_{\an}=2^n+2^{n-1}+\dots+2^m=2^{n+1}-2^m
$$
for some $m$ with
$3\leq m\leq n-1$ and $\dim(\phi)=2^n+2^{n-1}+\dots+2^m+2\iw_1$ where
$\iw_1=\iw_1(\phi)$ is the first Witt index of $\phi$.
Note that $\iw_1<2^{m-1}$ simply because $\dim(\phi)<2^{n+1}$.
Now it follows by \cite[conject. 0.1]{i_1}
or by Item 1 of Proposition \ref{known}
(take in account that the highest Witt index of $\phi$ is $2^{n-1}$)
that $\iw_1=2^{p-1}$ for some integer $p$ satisfying
$1\leq p\leq m-1$, and $\dim(\phi)=2^n+2^{n-1}+\dots+2^m+2^p$.
Since $\phi$ is a counter-example, $p$ is not $m-1$, so that
$1\leq p\leq m-2$ in fact.

Finally, \cite[conject. 0.1]{i_1}
(or Item 1 of Proposition \ref{known})
and the fact that
all
dimensions $\dim(\phi_E)_0<\dim(\phi)$ are allowed by Conjecture
\ref{conj},
allows one to determine all further higher Witt indices
of $\phi$ (compare with the proof of Corollary \ref{s.p.small}).
They are as follows (starting from $\iw_2$):
$2^{m-1}, 2^{m-2}, \dots, 2^{n-1}$,
meaning that the splitting pattern of $\phi$ consists of the
partial sums of the sum $2^n+2^{n-1}+\dots+2^m+2^p$.
Therefore the hypothesis of the preceding {\S} (and, in particular, the hypothesis of
Proposition \ref{mu} and Lemma \ref{mu'}) are satisfied.

To simplify the formulae which follow, like we did in the previous {\S},
we introduce the notation
$$
a=\iw_1(\phi)=2^{p-1};\;\;
b=\iw_2(\phi)=2^{m-1};\;\;
c=\iw_3(\phi)=2^m
$$
and
$$
d=\dim(X)/2=2^{n-1}+2^{n-2}+\dots+2^{m-1}+2^{p-1}-1\;.
$$

Let us consider the cycle $\mu\in\BCH(X^3)$ of Proposition \ref{mu}
as a correspondence from $\bar{X}$ to $\bar{X}^2$;
let us consider the cycle
$S^{2a}(\mu)\cdot(h^0\times h^0\times h^{b-1})\in\BCH(X^3)$ as a
correspondence from $\bar{X}^2$ to $\bar{X}$
(where $S^i$ stays for the degree $i$ component of the total Steenrod operation $S$).
Then we may take the composition of correspondences
$$
\mu\compose
\big(S^{2a}(\mu)\cdot(h^0\times h^0\times h^{b-1})\big)
\in\BCH(X^4)\;.
$$
Let us additionally consider
the morphism
$$
\delta\!:X^2=X_1\times X_2\to
X_1\times X_2\times X_3\times X_4=X^4\;,\;\;
x_1\times x_2\mapsto x_1\times x_2\times x_1\times x_2
$$
(all $X_i$ are copies of $X$)
given by the product of the diagonals
$X_1\to X_1\times X_3$ and $X_2\to X_2\times X_4$
(that is, $\delta$ is the diagonal morphism of $X^2$).

The following proposition
contradicts to Proposition \ref{known}
(note that $\dim(\xi)=2d+b-2a-1\geq2d+a$)
and proves therefore
Conjecture \ref{conj}.

\begin{prop}
\label{contradiction}
Let $\mu\in\BCH(X^3)$ be the cycle of Proposition \ref{mu}.
Then the decomposition of the cycle
$$
\xi=
\delta^*\Big(\mu\compose\big(S^{2a}(\mu)\cdot(h^0\times h^0\times h^{b-1})\big)\Big)
\in\BCH_{2d+b-2a-1}(X^2)
$$
contains the basis cycle $h^a\times l_{b-a-1}$.
\end{prop}

\begin{proof}
It is easy to see that each power of $h$ which is a factor of a basis element
involved in the decomposition of the cycle
$\mu\compose\big(S^{2a}(\mu)\cdot(h^0\times h^0\times h^{b-1})\big)$
is a multiple of $a$.
Therefore the same is true for the cycle $\xi$.

As before, we set $\mu_0=\mu-\mu'$.
We have:
\begin{multline*}
\xi=
\delta^*\Big(\mu_0\compose\big(S^{2a}(\mu_0)\cdot(h^0\times h^0\times h^{b-1})\big)\Big)
+\\
\delta^*\Big(\mu'\compose\big(S^{2a}(\mu')\cdot(h^0\times h^0\times h^{b-1})\big)\Big)
+\\
\delta^*\Big(\mu'\compose\big(S^{2a}(\mu_0)\cdot(h^0\times h^0\times h^{b-1})\big)\Big)
+\\
\delta^*\Big(\mu_0\compose\big(S^{2a}(\mu')\cdot(h^0\times h^0\times h^{b-1})\big)\Big)\;,
\end{multline*}
and we consider each of these four summands separately, one by one.

\bigskip
\noindent
{\bf First summand.}
First of all we compute $S^{2a}(\mu_0)$.
For every summand $h^0\times h^{(i-1)b+a}\times l_{ib+a-1}$ in the decomposition of
$\mu_0$, and for any $r\ne a$ with $1\leq r\leq2a$, we have:
$S^r(h^{(i-1)b+a})=S^r(l_{ib+a-1})=0$
(for $r=2a$ the relation $4a|b$, that is, $b\ne2a$, is used)
while
$S^a(h^{(i-1)b+a})=h^{(i-1)b+2a}$ and $S^a(l_{ib+a-1})=l_{ib-1}$.
Therefore
$$
S^{2a}(h^0\times h^{(i-1)b+a}\times l_{ib+a-1})=
h^0\times h^{(i-1)b+2a}\times l_{ib-1}
$$
and
$$
S^{2a}(\mu_0)=
\Sym\left(\Sum_{i=1}^{(d-a+1)/b}h^0\times h^{(i-1)b+2a}\times l_{ib-1}\right)\;.
$$
Now we can calculate the composition
$$
\mu_0\compose
\big(S^{2a}(\mu_0)\cdot(h^0\times h^0\times h^{b-1})\big)
$$
(see Lemma \ref{comp-corr}).
The basis cycles which appear in the decomposition of
$S^{2a}(\mu_0)\cdot(h^0\times h^0\times h^{b-1})$
have on the third factor place the following elements:
\begin{equation}
\tag{$*$}
h^{b-1},\;\; h^{ib+2a-1},\;\; l_{(i-1)b}\;.
\end{equation}
On the other hand, the basis elements which appear in the decomposition
of $\mu_0$ itself have on the {\em first} factor place the following:
\begin{equation}
\tag{$**$}
h^0,\;\; h^{(i-1)b+a},\;\; l_{ib+a-1}\;.
\end{equation}
It is straight-forward to see that the only pair of elements,
one from $(*)$, one from $(**)$, with the product $l_0$ is
$(l_0,h^0)$
(look at the indices modulo $2a$).
Therefore
\begin{multline*}
\mu_0\compose
\big(S^{2a}(\mu_0)\cdot(h^0\times h^0\times h^{b-1})\big)=
\\
\left(h^0\times
\Sym\left(\Sum_{i=1}^{(d-a+1)/b}h^{(i-1)b+a}\times l_{ib+a-1}\right)\right)
\compose
\big(\Sym(h^0\times h^{2a})\times l_0\big)=\\
\Sym(h^0\times h^{2a})\times
\Sym\left(\Sum_{i=1}^{(d-a+1)/b}h^{(i-1)b+a}\times l_{ib+a-1}\right)\;.
\end{multline*}
Applying $\delta^*$ to the cycle obtained, we get
\begin{multline*}
\Sym\left(\Sum_{i=1}^{(d-a+1)/b}h^{(i-1)b+a}\times l_{ib-a-1}+
h^{(i-1)b+3a}\times l_{ib+a-1}\right)=
h^a\times l_{b-a-1}+\dots\;.
\end{multline*}
It remains to show that the ``remaining part'' of $\xi$ does not contain the
basis cycle $h^a\times l_{b-a-1}$.

\bigskip
\noindent
{\bf Second summand.}
A basis cycle of the shape
$h^x\times?\times h^y\times?$
can be involved in the decomposition of
$$
\mu'\compose
\big(S^{2a}(\mu')\cdot(h^0\times h^0\times h^{b-1})\big)
$$
only if $x,y\geq a$.
In this case $\delta^*(h^x\times?\times h^y\times?)=h^{x+y}\times?$
with $x+y>a$,
therefore the cycle
$$
\delta^*\Big(\mu'\compose
\big(S^{2a}(\mu')\cdot(h^0\times h^0\times h^{b-1})\big)\Big)
$$
does not contain the basis element $h^a\times l_{b-a-1}$.

\bigskip
\noindent
{\bf Third summand.}
In order to check that the cycle
$$
\delta^*\Big(\mu'\compose
\big(S^{2a}(\mu_0)\cdot(h^0\times h^0\times h^{b-1})\big)\Big)
$$
does not contain the basis element $h^a\times l_{b-a-1}$,
it suffices to check that the number of basis elements in the decomposition
of the composition
$$
\mu'_2\compose
\big(S^{2a}(\mu_1)\cdot(h^0\times h^0\times h^{b-1})\big)
$$
is even
(note that we replaced $\mu'$ by $\mu'_2$, the notation being introduced in
Lemma \ref{mu'}, and we replaced $\mu_0$ by $\mu_1$,
the notation $\mu_1$ being introduced in the proof of Proposition \ref{mu} for
the sum of the basis elements contained in the decomposition of $\mu$ having
$h^0$ on the first factor place).
For this,
due to Lemma \ref{mu'},
it suffices to check that
every of the $(c-b)/a$ compositions
$$
t_{12}(\chi_j)\compose
\big(S^{2a}(\mu_1)\cdot(h^0\times h^0\times h^{b-1})\big)\;,
\;\;j\in[1,\;(c-b)/a]
$$
contains an even number of basis elements.
We show this by a straight-forward computation.
The point is that the number of summands in the decomposition of every $\chi_j$
is even
and either each or none of the summands ``produces''
a basis element
in the composition
(moreover, in the first case, precisely one basis element is produced by each summand
of the cycle $\chi_j$).

Let us do the computation.
The cycles $S^{2a}(\mu_1)\cdot(h^0\times h^0\times h^{b-1})$
and $\chi_j$ (for a fixed $j$)
are equal respectively to $h^0\times\alpha$ and to $h^a\times \beta$,
where
$$
\alpha=\Sum_{i=1}^{(d-a+1)/b}
h^{(i-1)b+2a}\times l_{(i-1)b}+l_{ib-1}\times h^{ib+2a-1}\;,
$$
while
$$
\beta=\Sum_{i=1}^{(d-b-a+1)/c}
h^{(i-1)c+b+ja}\times l_{(i-1)c+2b+(j+1)a-1}+
l_{ic+b-(j-2)a-1}\times h^{ic-(j-1)a}\;,
$$
and we just need to check that the composition $\beta\compose\alpha$ is
a sum of an even number of basis cycles.
There are two different cases depending on the value of $j$.
If the product $ja$ is not $0$ modulo $b$, then every product of every second factor of
the basis cycles appearing in the decomposition of $\alpha$ (namely,
$l_{(i-1)b}$ and $h^{ib+2a-1}$) by every first factor of
the basis cycles appearing in the decomposition of $\beta$
(namely, $h^{(i-1)c+b+ja}$ and $l_{ic+b-(j-2)a-1}$)
is different from $l_0$
(to see this, look at the indices modulo $b$).
Therefore, the composition of every basis cycle appearing
in the decomposition of $\alpha$ with every
basis cycles appearing in the decomposition of $\beta$
is $0$, and so, $\beta\compose\alpha=0$ in this case.

In the contrary case --- the case with $ja\equiv0\pmod{b}$ --- for every basis
cycle $y$ in the decomposition of $\beta$ there is precisely one basis cycle $x$ in the
decomposition of $\alpha$ such that $y\compose x\ne0$ (note that $y\compose x$
is a basis cycle in this case).
Since the number of basis cycles in the decomposition of $\beta$ is even
(equal to the integer
$(d-b-a+1)/c$ doubled), the composition $\beta\compose\alpha$ is the sum of an even number
of basis cycles.

\bigskip
\noindent
{\bf Fourth summand.}
We finish the proof of Proposition \ref{contradiction} considering the cycle
$$
\delta^*\Big(\mu_0\compose
\big(S^{2a}(\mu')\cdot(h^0\times h^0\times h^{b-1})\big)\Big)\;.
$$
We replace $\mu'$ by $\mu'_1$ and, furthermore, $\mu'_1$ by $\chi_j$ with some
$j\in[1,\;(c-b)/a]$.
Also, we replace $\mu_0$ by $\mu_2$.
We are going to show that the number of the basis elements
of the shape
$h^a\times?\times h^0\times?$
appearing in the decomposition of the composition
$$
\mu_2\compose
\big(S^{2a}(\chi_j)\cdot(h^0\times h^0\times h^{b-1})\big)
$$
is even (this will finish the proof of Proposition \ref{contradiction}).
We have $\chi_j=h^a\times\beta$ and $\mu_2=t_{12}(h^0\times\alpha)$
with
$$
\alpha=
\Sum_{i=1}^{(d-a+1)/b}h^{(i-1)b+a}\times l_{ib+a-1}+
l_{ib+a-1}\times h^{(i-1)b+a}\;,
$$

$$
\beta=
\Sum_{i=1}^{(d-b-a+1)/c}
h^{(i-1)c+b+ja}\times l_{(i-1)c+2b+(j+1)a-1}+
l_{ic+b-(j-2)a-1}\times h^{ic-(j-1)a}\;.
$$
Therefore the number we are looking for is the number of summands in the decomposition
of $\alpha\compose \big(S^{2a}(\beta)\cdot(h^0\times h^{b-1})\big)$.

We can compute the Steenrod operation $S^{2a}$ on the summands of the decomposition
of $\beta$.
The formula depends on the value of $j$ modulo $4$ because of the rules
(here $S^{\leq2a}$ stays for $\Sum_{k\leq2a}S^k$):
$$
S^{\leq2a}(h^{ia})=
\begin{cases}
h^{ia} &\text{if $i\equiv0\pmod{4}$;}\\
h^{ia}+h^{(i+1)a}  &\text{if $i\equiv1\pmod{4}$;}\\
h^{ia}+\hspace{4em}h^{(i+2)a} &\text{if $i\equiv2\pmod{4}$;}\\
h^{ia}+h^{(i+1)a}+h^{(i+2)a} &\text{if $i\equiv3\pmod{4}$,}
\end{cases}
$$
while
(here recall that $S(l_i)=l_i\cdot(1+h)^{2d-i+1}$, $d+1\equiv a\pmod{b}$,
and $4a$ divides $b$):
$$
S^{\leq2a}(l_{ia-1})=
\begin{cases}
l_{ia-1}+\hspace{4.5em}l_{(i-2)a-1} &\text{if $i\equiv0\pmod{4}$;}\\
l_{ia-1}+l_{(i-1)a-1}  &\text{if $i\equiv1\pmod{4}$;}\\
l_{ia-1} &\text{if $i\equiv2\pmod{4}$;}\\
l_{ia-1}+l_{(i-1)a-1}+l_{(i-2)a-1} &\text{if $i\equiv3\pmod{4}$.}
\end{cases}
$$

Assume that $j\equiv0\pmod{4}$ or $j\equiv1\pmod{4}$.
Then, applying the above formulae, we get that
$S^{2a}(\beta)=0$,
and there is nothing more to prove in this case.

Now we assume that $j\equiv2\pmod{4}$ or $j\equiv3\pmod{4}$.
Then
$$
S^{2a}(\beta)\cdot(h^0\times h^{b-1})=\beta_1+\beta_2\;,
$$
where
\begin{multline*}
\beta_1=\beta\cdot(h^{2a}\times h^{b-1})=\\
\Sum_{i=1}^{(d-b-a+1)/c}
h^{(i-1)c+b+(j+2)a}\times l_{(i-1)c+b+(j+1)a}+
l_{ic+b-ja-1}\times h^{ic+b-(j-1)a-1}\;,
\end{multline*}
while
\begin{multline*}
\beta_2=\beta\cdot(h^0\times h^{b+2a-1})=\\
\Sum_{i=1}^{(d-b-a+1)/c}
h^{(i-1)c+b+ja}\times l_{(i-1)c+b+(j-1)a}+
l_{ic+b-(j-2)a-1}\times h^{ic+b-(j-3)a-1}\;.
\end{multline*}

If $j\equiv3\pmod{4}$, then the compositions $\alpha\compose\beta_1$ and
$\alpha\compose \beta_2$ are $0$
because so are the compositions of any basis cycles included in $\alpha$
with any basis cycle included in $\beta_1$ or $\beta_2$
(look at the indices modulo $2a$).
If $j\equiv2\pmod{4}$, then $\alpha\compose\beta_1=0$ too and by the same reason
(look at the indices modulo $4a$).
Finally, assume that $j\equiv2\pmod{4}$ and
consider the composition
$\alpha\compose\beta_2$.
For every basis
cycle $y$ in the decomposition of $\beta_2$ there is precisely one basis cycle $x$ in the
decomposition of $\alpha$ such that $x\compose y\ne0$ (note that $x\compose y$
is a basis cycle in this case).
Since the number of basis cycles in the decomposition of $\beta_2$ is even
(equal to the integer
$(d-b-a+1)/c$ doubled), the composition $\alpha\compose\beta_2$ is the sum of an even number
of basis cycles.
\end{proof}

Conjecture \ref{conj} is proved.
The following supplement is now easy to get:

\begin{cor}
\label{s.p.small}
Let $\phi$ be a small quadratic form with
$\dim(\phi)=2^{n+1}-2^m$, $m\in[1,\; n+1]$.
Then the splitting pattern
$\{\dim(\phi_E)_{\an}|\;\text{$E/F$ is a field extension}\}$
of the form $\phi$ coincides with the set $\{2^{n+1}-2^i\}_{i=m}^{n+1}$
(in particular, the height of $\phi$ is equal to $n+1-m$).
\end{cor}

\begin{proof}
By Conjecture \ref{conj} proved above,
$\dim(\phi_{F(X)})_{\an}=2^{n+1}-2^r$ for some $r\in[m+1,\; n+1]$.
But by \cite[conject. 0.1]{i_1} (or by Item 1 of Proposition \ref{known}
taking in account that the highest Witt index of $\phi$ is
$2^{n-1}$), it follows that only the value $r=m+1$ is possible.
Proceeding this way (with the form $(\phi_{F(X)})_{\an}$ and so on), we get the result.
\end{proof}

\section
{Possible dimensions}
\label{Possible dimensions}

Let us recall some standard  notation concerning quadratic forms:
one writes $\<a_1,\dots,a_n\>$, where $a_1,\dots,a_n\in F$, for the
quadratic form
$$
F^n\to F\;,\;\; (x_1,\dots,x_n)\mapsto a_1x_1^2+\dots+a_nx_n^2\;;
$$
$\ll a_1,\dots,a_n\rr$ for the $n$-fold Pfister form
$$
\<1,-a_1\>\otimes\dots\otimes\<1,-a_n\>\;,
$$
and $\ll a_1,\dots,a_n\rr'$ for the pure subform of the above Pfister
form (see \cite[def. 1.1 of ch. 4]{Scharlau}).

The following result provides, in particular, examples for all dimensions
which are not prohibited by Conjecture \ref{conj}.

\begin{thm}[A. Vishik]
\label{examples}
Take any integers $n\geq1$ and $m\geq2$.
Let $k$ be a field (of $\Char(k)\ne2$),
$t_i$, $t_{ij}$ ($i=1,\dots,m$, $j=1,\dots,n$) variables, and
$$
F=k(t_i,t_{ij})_{1\leq i\leq m,\; 1\leq j\leq n}
$$
the field of rational functions in all these variables.

The splitting pattern of the quadratic form
$$
\phi=t_1\cdot\ll t_{11},\dots,t_{1n}\rr\bot\dots\bot
\;t_m\cdot\ll t_{m1},\dots,t_{mn}\rr
$$
over $F$ is
$$
\{2^{n+1}-2^i|\;i=n+1,n,\dots,1\}\cup\big(2\Z\cap[2^{n+1},\;m\cdot2^n]\big)\;.
$$
\end{thm}

\begin{proof}
First of all, it is easy to see that all the integers $2^{n+1}-2^i$ are in the
splitting pattern of $\phi$.
Indeed, the anisotropic part of $\phi$ over the field $E$ obtained from $F$
by adjoining the square roots of $t_{31},t_{41}\dots,t_{m1}$, of $t_1$  and of
$-t_2$, is isomorphic to the (generalized Albert) form
$$\ll t_{11},\dots, t_{1n}\rr'\bot-
\ll t_{21},\dots,t_{2n}\rr'
$$
(the primes $'$ mean the pure subforms of the Pfister forms)
of dimension $2^{n+1}-2$;  the splitting pattern of this form is
$\{2^{n+1}-2^i\}$ because this set is the splitting pattern of {\em any}
anisotropic
$(2^{n+1}-2)$-dimensional quadratic form whose class lies in
$I^n$ (Corollary \ref{s.p.small}).

Now let us assume that some (at least one) even integers of the interval
$[2^{n+1},\;m\cdot2^n]$ are {\em not} in the splitting pattern of $\phi$.
Among all such integers we take the smallest one and call it $a$;
let $b$ be the biggest integer smaller that $a$ and lying in the splitting
pattern;
let $c$ be the smallest integer greater that $a$ and lying in the splitting
pattern.
Let $E$ be the field of the generic splitting tower of $\phi$ such that
$\dim(\psi)=c$ for $\psi=(\phi_E)_{\an}$.
Let $Y$ be the projective quadric given by the quadratic form $\psi$.
Let $\pi\in\BCH(Y^2)$ be the cycle of the set $\Pi$ of Theorem \ref{kruto!}
with $f(\pi)=1$.
We claim that $\pi=\Sym(h^0\times l_{\iw_1-1})$ for $\iw_1=\iw_1(Y)$.
Indeed, since $\iw_1=(c-b)/2>1$ and $\iw_q(Y)=1$ for all
$q\in S=S(Y)=\{1,2,\dots,\Hight\}$
($\Hight$ is the height of $\psi$) such that
$\dim(\psi_q)\in[2^{n+1}-2,\;b]$, the cycle
$\pi$ does not contain $h^{\iw_1+\dots+\iw_{q-1}}\times
l_{\iw_1+\dots+\iw_{q-1}+\iw_1-1}$ for such $q$
(Lemma \ref{dyrki v verhnem cykle}), and also
$\pi\not\ni h^i\times l_{i+\iw_1-1}$ for any $i\in\{1,2,\dots,\iw_1-1\}$
(Lemma \ref{dyrki v verhnem cykle} as well).
Finally, for the integer $q\in S$ such that
$\dim(\psi_q)=2^{n+1}-2$,
the cycle $\pr^2(\pi_{E_q})\in\BCH(Y_q^2)$
(the homomorphism $\pr^2\!:\BCH(Y^2_{E_q})\to\BCH(Y_q^2)$ is defined in
Corollary \ref{old isotropic}) has the dimension
$$
\dim(Y_q)+\iw_1-1\geq\dim(Y_q)+1=\dim(Y_q)+\iw_1(Y_q)
$$
and therefore is $0$ by Item \ref{item 3} of Proposition \ref{known}.

We have shown that $\pi=\Sym(h^0\times l_{\iw_1-1})$.
By Item \ref{binary size} of Proposition \ref{list}, it
follows that the integer $\dim(Y)-\iw_1+1$ is a power of $2$, say
$2^p$.
Since
$$
\dim(Y)-\iw_1+1=(c-2)-(c-b)/2+1=(b+c)/2-1\;,
$$
the integer $2^p$ sits inside of the open interval
$(b,\;c)$; therefore, satisfying
$2^{n+1}\leq2^p<m\cdot2^n$, the integer $2^p$
is {\em not} in the splitting pattern of the
quadratic form $\phi$.
But all the integers $\leq m\cdot2^n$ divisible by $2^n$ are evidently {\em in} the
splitting pattern of $\phi$.
The contradiction obtained proves Theorem \ref{examples}.
\end{proof}

\begin{rem}
\label{easy examples}
Of course, the dimensions $2^{n+1}-2^i$ can be realized more directly by
the tensor products of Pfister forms and generalized Albert forms
($u_.,v_.,w_.$ are variables):
$$
\ll u_1,\dots u_{i-1}\rr\otimes\big(\ll v_1,\dots,v_{n+1-i}\rr'
\bot-\ll w_1,\dots,w_{n+1-i}\rr'\big)\;.
$$
\end{rem}

\end{document}